\documentclass[a4paper,11pt,reqno]{amsart}
\usepackage{amsthm}
\usepackage{amsmath}
\usepackage{enumerate}
\usepackage{amsfonts}
\usepackage{amssymb}
\usepackage{fullpage}
\usepackage{amsmath,amscd}
\usepackage{stmaryrd}
\usepackage{graphicx}
\usepackage{tikz-cd}
\usepackage{array}
\usepackage{amsmath,amssymb,amsfonts,dsfont}
\usepackage[utf8]{inputenc}
\usepackage[english]{babel}
\usepackage[T1]{fontenc}
\usepackage{float}
\usepackage{enumerate}
\usepackage{pdfpages}
\usepackage[a4paper, margin = 3cm, bottom = 3cm]{geometry}
\usepackage{ifpdf}
\usepackage{marginnote}
\usepackage{tikz}	
\usepackage{tikz-cd}
\graphicspath{{images/}}
\usepackage{esvect}
\usepackage{yhmath} 
\usepackage{placeins}
\usepackage{hyperref} 

\definecolor{arnaudcolor}{rgb}{0.5 0 1}

\usepackage[normalem]{ulem}

\usetikzlibrary{calc}
\usetikzlibrary{decorations.pathreplacing,decorations.markings,decorations.pathmorphing}
\usetikzlibrary{positioning,arrows,patterns}
\usetikzlibrary{cd}
\usetikzlibrary{intersections}
\usetikzlibrary{arrows}

\hypersetup{pdfborder=0 0 0,
	    colorlinks=true,
	    citecolor=black,
	    linkcolor=blue,
	    urlcolor=red,
	    pdfauthor={Guillaume Tahar}
	   }
\renewcommand{\Re}{\operatorname{Re}}
\renewcommand{\Im}{\operatorname{Im}}

\newcommand{\CC}{{\mathbb{C}}}
\newcommand{\C}{{\mathbb{C}}}
\newcommand{\RR}{{\mathbb{R}}}
\newcommand{\R}{{\mathbb{R}}}
\newcommand{\Id}{{\on{Id}}}
\newcommand{\Ric}{\mathrm{Ric}}
\newcommand{\wt}{\widetilde}
\newcommand{\eps}{\epsilon}
\newcommand{\ora}{\vec}

\newcommand{\on}[1]{\operatorname{#1}}
\newcommand{\cal}[1]{\mathcal{#1}}

\newtheorem{thm}{Theorem}[section]
\newtheorem{cor}[thm]{Corollary}
\newtheorem{prop}[thm]{Proposition}
\newtheorem{lem}[thm]{Lemma}

\theoremstyle{definition}
\newtheorem{defn}[thm]{Definition}
\theoremstyle{remark}
\newtheorem{rem}[thm]{Remark}
\theoremstyle{definition}

\theoremstyle{definition}
\newtheorem{ex}[thm]{Example}
\theoremstyle{definition}

\numberwithin{equation}{section}
\title{Hybrid connections on Hessian manifolds}

\author{Arnaud Chéritat}
\address[Arnaud Chéritat]{Univ Toulouse, INUC ou UT2J, INSA, CNRS, IMT, Toulouse, France}
\email{arnaud.cheritat@math.univ-toulouse.fr}

\author{Guillaume Tahar}
\address[Guillaume Tahar]{Department of Mathematics with Computer Science, Guangdong Technion - Israel Institute of Technology (GTIIT), Shantou, Guangdong, China}
\email{guillaume.tahar@gtiit.edu.cn}

\date{\today}

\begin{document}

\begin{abstract}
We introduce a new class of affine connections on Hessian manifolds, called hybrid connections, characterized by the compatibility between their projective geometry with the underlying affine structure on the one hand and their infinitesimal holonomy with the Hessian metric on the other hand. 
In this paper, we investigate the properties of hybrid connections and prove that, on a given Hessian manifold, they are completely determined by the choice of a Hessian potential for the metric.
\par
In the special case of pseudo-Euclidean manifolds, we identify canonical models and construct, in particular, a natural connection on the open unit ball that combines features of the Cayley–Klein and Poincaré models of hyperbolic geometry. We also prove the existence and uniqueness (up to scaling) of a pseudo-Riemannian metric $h$ such that the geodesics of $\nabla$ admit parameterizations of constant speed with respect to $h$, which we call the isochrone metric.
\end{abstract}

\maketitle

\let\oldtocsubsection=\tocsubsection
\renewcommand{\tocsubsection}[2]{\hspace{1em}\oldtocsubsection{#1}{#2}}
\setcounter{tocdepth}{2}
\tableofcontents

\section{Introduction}\label{sec:intro}

In this paper, we investigate a form of compatibility between projective flatness and infinitesimal conformal flatness, for affine connections.

In the whole article, \emph{metric} stands for \emph{pseudo-Riemannian} metric.

\subsection{In search of the best connection}

Let us first present a motivating example, before defining in Section~\ref{sub:hc} a general setting in which it fits. 

\medskip

Planar charts of non-Euclidean surfaces are always unsatisfactory in one way or another.
One well-known reason is that the Gaussian curvature prevents the chart to map the metric of the surface to the canonical Euclidean metric in the chart.

In some case, we can recover a partial match between the properties of the two metrics. For example, there are two famous representations of the hyperbolic plane as the unit disk $D\subset\C$: the Poincaré model, whose metric is, up to multiplication by a constant, the unique one which is invariant by the group of complex homographies preserving $D$ and the Beltrami-Klein model, whose model is, up to multiplication by a constant,  the unique one that is invariant by the group of projective transformations of $\mathbb{RP}^2$ preserving $D$.
See Table~\ref{tab:Comparison} for analytic expressions of their respective metrics.
The first one is conformal (angles between curves are identical for the hyperbolic metric and for the Euclidean one) but hyperbolic geodesics appear as curved lines.
In the second model, the geodesics of hyperbolic geometry are represented as the chords of the disk but the Euclidean angle at their intersection does not faithfully represent the true one.
In a similar way, gnomonic projection of the sphere sends great circles to straight lines, but it is not a conformal representation either, while no conformal charts can send all geodesics to segments of straight lines.

As a matter of fact, on a connected open subset of the Euclidean plane, the only metrics that are conformal to the Euclidean metric and whose geodesics are straight lines are the Euclidean metric itself and its multiplication by a non-zero constant.

\smallskip

By focusing on connections rather that metrics, one gets more flexibility.
To a metric $g$ is associated the Levi-Civita connection, which is the unique symmetric affine connection that preserves $g$.
Its geodesics are the same as the geodesics of $g$, and one can recover $g$, up to multiplication by a constant, from the connection.

Let now take for $g$ the Euclidean metric in the plane.
As we saw, the only symmetric connection preserving $g$ is the canonical one.
Let us relax the condition that the connection $\nabla$ should $g$ to the following weaker pair of requirements: (1) the first-order infinitesimal holonomies of $\nabla$ are infinitesimal isometries of $g$, (2) the geodesics of $\nabla$ are straight lines of the plane.
Surprisingly, we obtain a family of non-trivial connections with remarkable properties.
On open subsets of the Euclidean plane, they are given, up to a translation, by the formula $\Gamma(u,v)=a(v)u+a(u)v$ for their Christoffel symbol (as is every symmetric connection whose geodesics are straight lines), where $a$ is the differential form with the following simple expression:
\begin{equation}\label{eq:Example1}
a(u) = - \frac{2x\,dx+2y\,dy}{x^{2}+y^{2}+\lambda}.
\end{equation}
Here $x,y$ are the Cartesian coordinates of $u=(x,y)$ and $\lambda$ is an arbitrary real constant.
We will see in this article how such a formula appears.
Specializing to $\lambda = -1$, we get a new model of geometry on the unit disk, which is discussed in details in Section~\ref{sec:Tale}.
Interestingly, it is related to the Riemannian metric $h = |dz|^2/(1-|z|^2)^4$ with $z=x+iy$: the geodesics of the connection have constant speed with respect to $h$.
We say that the connection $\nabla$ and the metric $h$ are \emph{isochrone}.
Note that $h$ is not the Poincaré metric, the latter having expression $|dz|^2/(1-|z|^2)^2$.
Of course, $h$ and $\nabla$ have different geodesics: those of $\nabla$ are straight lines while those of $h$ are shown on the left part of Figure~\ref{fig:geod:ex}.

\subsection{Hybrid connections}\label{sub:hc}

In Equation~\eqref{eq:Example1}, we recognize the expression $-\frac{d\Psi}{\Psi}$ where $\Psi$ is a Hessian potential of the standard Euclidean metric $dx^{2}+dy^{2}$. The natural framework in which the examples mentioned above are generalized is that of Hessian manifolds, which combine both an affine structure and a metric, and whose definition we recall below.
\par
Hessian metrics are a real variable analog of K\"{a}hler metrics for complex variable, and arise in the study of the real Monge-Ampère equation (see \cite{cheng1982}). More recently, Hessian metrics appeared in information geometry as the Fisher-Rao metric of parameter spaces of some statistical models (see \cite{Shima1981,Shima1997} for background). 

An \emph{affine manifold} is a pair $(M,D)$ with $M$ a manifold and $D$ a flat affine connection. It has \emph{affine charts} in which $D$ is trivial (equal to the canonical connection of $\R^d$).
It also carries a \emph{Hessian operator}, $H=D^2$, which has a very simple expression in affine charts:
\[ H f = \sum_{i,j} \frac{\partial^2 f}{\partial x_i\partial x_j} dx_i\,dx_j
,\]
to be understood as a symmetric bilinear form on the tangent space (a symmetric $(0,2)$-tensor).

\begin{defn}\label{defn:Hessian}
A Hessian manifold $(M,D,g)$ is a triple where:
\begin{itemize}
    \item $M$ is a connected real manifold of dimension $d \geq 2$;
    \item $D$ is a flat affine connection on $M$;
    \item $g$ is a pseudo-Riemannian metric of arbitrary signature $(p,q)$ with $p+q=d$ that has a local expression of the form $g=H \psi$ where $\psi$ is a (locally defined) smooth function called a \textit{Hessian potential}.
\end{itemize}
Note that since $M$ is connected, the non-degeneracy condition on $g$ implies that the signature $(p,q)$ is constant on $M$.
\end{defn}

In this paper (Definition~\ref{defn:Hybrid}), we introduce the notion of \textit{hybrid connection} where the unparameterized geodesics are determined by the affine structure induced by $D$ while the infinitesimal holonomy is governed by the pseudo-Riemannian metric $g$. 
\par
This notion, which emerged from an exploration aimed at finding generalized notions of flatness, fits into a broader theme in differential geometry: the search for distinguished classes of connections singled out by natural geometric constraints, as illustrated by Yang–Mills or Einstein connections.
\par
A related motivation comes from the growing study of flat geometric structures (such as translation, dilation, or affine structures) and their moduli spaces. These objects share a similar dual feature: geodesics analogous to those of an affine space, together with holonomy constraints that impose a rigid geometric structure.
\par
The notion of a hybrid connection is remarkable in that it gives rise (see Theorem~\ref{thm:MAIN}) to a finite-dimensional deformation space, identified with the space of globally defined Hessian potentials of the Hessian metric. These connections admit explicit expressions, and their parallel transports can be computed precisely, which we carry out, at least locally, in the case of pseudo-Euclidean spaces in Sections~\ref{sec:LocalModels} and~\ref{sec:Tale}.
\par
Before introducing formally the notion of a hybrid connection, we first need two notions: incompressible connections and infinitesimally conformal connections with respect to pseudo-Riemannian metrics.

\begin{defn}\label{defn:Incompressible}
A connection is called \emph{incompressible} whenever there exists, near every point, a non-vanishing volume form that is locally preserved by the parallel transport.
\end{defn}

\begin{defn}\label{defn:InfinitesimallyConformal}
On an affine manifold $(M,D)$ endowed with a pseudo-Riemannian metric $g$, a connection $\nabla$ is called \emph{infinitesimally conformal} whenever in any affine chart, the holonomy of $\nabla$ around an arbitrarily small planar rectangle of area $\mathcal{A}$ is of the form $\Id+\mathcal{A}\cdot \mathcal{M}  +o(\mathcal{A})$ where $\mathcal{M}$ is an infinitesimal similarity of $g$ (see Section~\ref{sub:InfinitesimalHolonomy} for background).
\end{defn}

The notion of a hybrid connection combines the properties defined above in Definitions~\ref{defn:Incompressible} and~\ref{defn:InfinitesimallyConformal} with the additional condition that the unparameterized geodesics of the connection are straight lines in affine charts (and a non-degeneracy hypothesis on the curvature tensor to avoid pathological objects).

\begin{defn}\label{defn:Hybrid}
On a Hessian manifold $(M,D,g)$, a \textit{hybrid connection} $\nabla$ is a symmetric affine connection that satisfies the following properties:
\begin{enumerate}
    \item unparameterized geodesics of $\nabla$ and $D$ coincide ($\nabla$ is projectively equivalent to $D$);
    \item $\nabla$ is incompressible (see Definition~\ref{defn:Incompressible});
    \item $\nabla$ is infinitesimally conformal (see Definition~\ref{defn:InfinitesimallyConformal}); 
    \item the Riemann curvature tensor $R$ of $\nabla$ nowhere vanishes.
\end{enumerate}
\end{defn}

\begin{rem}\label{rem:equiv}
Point~(1) is equivalent to the geodesics of $\nabla$ being straight lines in affine charts of $M$. Points~(2) and~(3) can be formulated as: $\nabla$ is infinitesimally isometric, i.e.\ $\mathcal{M}$ is an infinitesimal isometry. 
\par
Indeed, if~(2) holds then the holonomies are trace-free, and a trace-free infinitesimally conformal endomorphism is an infinitesimal isometry.
Conversely if the hybrid connection $\nabla$ is infintesimally isometric, then consider the connection $\nabla^{\Lambda}$ that it induces on the determinant line bundle $\Lambda^{n}TM$. If infinitesimal holonomies 
are infinitesimal isometries, they are in particular trace-free so the infinitesimal holonomy of $\nabla^{\Lambda}$ vanishes which is equivalent to (2) because $\nabla^{\Lambda}$ is a line bundle ($\on{GL}_1(\R)$ is commutative).
\end{rem}

The main result of this paper is a description of hybrid connections by their \emph{hybrid potential}, through a formula that generalizes Equation~\ref{eq:Example1}.


\begin{defn}\label{defn:LogarithmicPotential}
Given an affine manifold $(M,D)$, we say that a function $\psi$ is a \textit{logarithmic potential} of a connection $\nabla$ if it satisfies $\nabla=D-(\frac{d\psi}{\psi}\otimes \operatorname{id} + \operatorname{id}\otimes\frac{d\psi}{\psi})$, i.e.
$$
\nabla_{u} v = \partial_{u} v -\frac{d\psi (u)v + d\psi(v) u}{\psi(x)}
$$
for any vector $u$ in the tangent space $T_{m}M$ of $M$ at some point $m \in M$, and for any vector field $v$ defined in a neighborhood of $m$,
where $d\psi$ in the expression is taken at the point $m$.
\end{defn}

The next theorem, proved in Section~\ref{sub:charHybrid},
is a classification of hybrid connections.
We show that each hybrid connection has for logarithmic potential one of the local Hessian potentials of the underlying Hessian manifold (they are defined up to the addition of an affine function). Actually, the fact that these local potentials simultaneously serve as logarithmic potentials for the connection and as Hessian potentials for the metric, both globally defined, forces them to be globally defined.

\begin{thm}\label{thm:MAIN}
On a Hessian manifold $(M,D,g)$, for any globally defined \textit{hybrid connection} $\nabla$, there exists a globally defined smooth function $\psi :M \to \R^*$ such that $H(\psi)=g$ and such that $\psi$ is a logarithmic potential of $\nabla$.
\par
Conversely, any connection of this form is a hybrid connection.
\end{thm}

Remarkably, $\psi$ is at the same time what we call a \textit{logarithmic potential} for $\nabla$ (defined up to a scaling, see Proposition~\ref{prop:Ricci}) and a \textit{Hessian potential} for $g$ (defined up to the addition of an affine function with respect to the affine connection $D$).
As a consequence, $\psi$ is uniquely defined.

\begin{defn}\label{def:hybridPotential}
We call this unique $\psi$ the \emph{hybrid potential} of $\nabla$.
\end{defn}

\begin{rem}
One could ask the more general question, given a (pseudo-)Riemannian metric $g$, whether or not there exists a symmetric connection whose infinitesimal holonomies are infinitesimal isometries of $g$ and which is projectively trivial (there is a chart where the geodesics of the connection are straight lines), however, this question is outside the scope of the present paper: we restrict ourselves to the case where $g$ is expressible as a Hessian metric in an affine chart and where the geodesics are straight lines in this chart.
As a benefit, we get this striking fact that the same function $\psi$ serves as a Hessian potential and a logarithmic one.
\end{rem}

\subsection{Pseudo-Euclidean manifolds}

\begin{defn}
A Hessian manifold $(M,D,g)$ is a \textit{pseudo-Euclidean manifold} if the pseudo-Riemannian metric $g$ is constant in the affine charts of $D$.
\end{defn}

This condition is invariant under affine changes of charts of $D$.

\smallskip

According to Theorem~\ref{thm:MAIN}, in the case of a pseudo-Euclidean manifold of signature $(p,d-p)$, the hybrid potential for a hybrid connection $\nabla$ can locally be written as
\[ \psi = \lambda + \frac{1}{2}\sum\limits_{i\leq p} {x_{i}}^{2} -\frac{1}{2}\sum\limits_{j > p} {x_{j}}^{2}
\] for some affine coordinates $x_{1},\dots,x_{d}$ and some constant $\lambda \in \mathbb{R}$.
\par
Although geodesics of $\nabla$ are usually not geodesics for the metric $g$, we prove that every geodesic of $\nabla$ has constant speed with respect to the pseudo-Riemannian metric $\frac{1}{\psi^{4}}\cdot g$, called the \textit{isochrone metric} associated with $(M,D,g)$ and $\nabla$ (see Section~\ref{sub:isochrone} for details).

\smallskip

A \emph{conformal class} is an equivalence class of pseudo-Riemannian metric $g$ up to multiplication by a positive function $f$, i.e.\ $g\sim f g$.
Remarkably, the Ricci tensor of the hybrid connection $\nabla$ (up to sign normalization) also belongs to the same conformal class as the pseudo-Euclidean metric $g$ and the isochrone metric $\frac{1}{\psi^{4}}\cdot g$.

\begin{thm}\label{thm:ThreeMetrics}
Given a hybrid connection $\nabla$ defined by a hybrid potential $\psi$ on a pseudo-Euclidean manifold $(M,D,g)$, the following three pseudo-Riemannian metrics on $\psi^{-1}(\mathbb{R}^{\ast})$ belong to the same conformal class:
\begin{itemize}
    \item the Hessian pseudo-metric $g$;
    \item $\frac{(d-1)}{|\psi|}g$, which coincides with the Ricci tensor $\Ric(\nabla)$ of $\nabla$ at points of $M$ where $\psi>0$ and with its opposite at points of $M$ where $\psi<0$;
    \item the \textit{isochrone metric} $\frac{1}{\psi^{4}}\cdot g$.
\end{itemize}
Moreover, $\nabla$ preserves the volume form $\frac{\mathrm{Vol}}{\psi^{d+1}}$ where $\mathrm{Vol}$ is the canonical volume form in any affine chart with respect to $D$.
\end{thm}

This theorem is proved in Propositions~\ref{prop:Ricci} and~\ref{prop:metr}.
An important part of this article is dedicated to the systematic description of hybrid connections defined on subsets of Euclidean spaces (which are a special case of Hessian manifolds). The fact that the metric is preserved by the infinitesimal holonomy of a hybrid connection does not prevent it from having a holonomy group much larger than the isometry group alone (see Section~\ref{sub:LocHol}): though this is a generic phenomenon for connections, those examples may serve as a remarkable pedagogical illustration.
\par
In the case of the Euclidean ball, it appears that exactly one hybrid connection is \textit{geodesically complete} (see Definition~\ref{def:geo:complete}).

\begin{thm}\label{thm:Uniqueness}
For any $d \geq 2$, there exists a unique \textit{geodesically complete} hybrid connection $\nabla$ on the Hessian manifold $(\mathcal{B}^{d},D,g)$ where $\mathcal{B}^{d}$ is the $d$-dimensional unit ball, $D$ is the standard affine connection and $g$ is the standard Euclidean metric.
Its hybrid potential is $\psi:(x_{1},\dots,x_{d}) \mapsto \frac{1}{2}\left(|x|^{2}-1\right)$.
\end{thm}

{ 
\def\arraystretch{1.5}
\begin{table}
    \centering
    \begin{tabular}{|c||c|c|c|c|}
    \hline
Connection  & Flat & Cayley-Klein & Poincaré & Hybrid \\
\hline
\hline
Preserved metric   
& $dx^{2}+dy^{2}$ & $\frac{dx^{2}+dy^{2}}{1-x^{2}-y^{2}}+\frac{(xdx+ydy)^{2}}{(1-x^{2}-y^{2})^{2}}$ 
& $\frac{4(dx^{2}+dy^{2})}{(1-x^{2}-y^{2})^{2}}$ 
& $\times$ \\
\hline
Preserved volume form   
& $dx \wedge dy$   
& $\frac{dx \wedge dy}{(1-x^{2}-y^{2})^{3/2}}$  
& $\frac{4\, dx \wedge dy}{(1-x^{2}-y^{2})^{2}}$  
& $\frac{dx \wedge dy}{(1-x^{2}-y^{2})^{3}}$ \\
\hline
Curvature form   
&  0  
& $\frac{-dx \wedge dy}{(1-x^{2}-y^{2})^{3/2}}$  
& $\frac{-4\,dx \wedge dy}{(1-x^{2}-y^{2})^{2}}$  
& $\frac{-2\, dx\wedge dy}{1-x^2-y^2}$ \\
\hline
Geodesically complete & $\times$ & \checkmark & \checkmark & \checkmark \\
\hline
Infinitesimally conformal & \checkmark & $\times$ & \checkmark & \checkmark \\
\hline
Straight geodesics & \checkmark & \checkmark & $\times$ & \checkmark \\
\hline
    \end{tabular}
    \vskip1ex
    \caption{Comparison of natural connections for the unit open disk. For readability, we adopt coordinates $(x,y)$ instead of $(x_1,x_2)$.}
    \label{tab:Comparison}
\end{table}
}

This theorem is proved in \ref{sub:geo:complete}.
The case $d=2$, gives us the unique geodesically complete hybrid connection $\nabla$ on the unit disk of the Euclidean plane.
In Table~\ref{tab:Comparison}, we compare $\nabla$ with the trivial connection and with the Levi-Civita connections of the usual geometric models of the hyperbolic disk:
\begin{itemize}
    \item The Cayley-Klein model (where geodesics are chords);
    \item The Beltrami-Poincaré model (where the Levi-Civita connection preserves the conformal structure of the disk).
\end{itemize}
Recall that the metric of the hyperbolic plane is geodesically complete. 
The geodesics of the Cayley-Klein model are the usual straight lines and thus coincide with (unparameterized) geodesics of $\nabla$. However, the isochrone metric of $\nabla$ does not coincide with the Cayley-Klein hyperbolic metric nor the Poincaré metric. Remarkably, in each of the four models presented in Table~\ref{tab:Comparison}, there is a preserved volume form $\frac{dx \wedge dy}{(1-x^{2}-y^{2})^{k}}$ where $k$ can be $0$, $\frac{3}{2}$, $2$ or $3$ depending on the model.

\subsection{Organization of the paper}

\begin{itemize}
    \item In Section~\ref{sec:AffCon}, we provide background on affine connections.
    \item In Section~\ref{sec:FlatAffine}, we give the expression, due to Weyl, in a manifold $M$ endowed with a flat affine connection $D$, of a connection $\nabla$ that is projectively equivalent to $D$ (i.e.\ the unparameterized geodesics of $\nabla$ and $D$ coincide). In this class, incompressible connections (connections locally preserving a volume form) are determined by the logarithmic derivative of a potential.
    \item In Section~\ref{sec:HessianMetrics}, we provide background on Hessian operators of real affine structures, on Hessian metrics and give explicit formulas for infinitesimal isometries of Hessian metrics in terms of Hessian potentials, proving in particular Theorem~\ref{thm:MAIN}.
    \item In Section~\ref{sec:LocalModels}, we define local models for hybrid connections in pseudo-Euclidean manifolds. We also provide explicit computations for the parallel transport in some subsets of the Euclidean and Minkowski planes (two-dimensional canonical models).    
    \item In Section~\ref{sec:Tale}, we present the most remarkable feature of hybrid connections on pseudo-Euclidean manifolds. In the same conformal class, we find the Hessian metric of the underlying Hessian structure, the Ricci tensor (or its opposite) of the hybrid connection $\nabla$, and an isochrone metric (with respect to which geodesics of $\nabla$ have constant speed). In particular, we prove Theorems~\ref{thm:ThreeMetrics} and~\ref{thm:Uniqueness}.
\end{itemize}

\paragraph{\bf Acknowledgements.} 
The authors thank Gabriel Khan and the anonymous referee for valuable remarks.

\paragraph{\bf Funding statement.}
Both authors are supported by the French National Research Agency under the project TIGerS (ANR-24-CE40-3604).

\section{Reminder on affine connections}\label{sec:AffCon}

In this paper, we study \textit{affine connections} on a smooth manifold $M$ (locally diffeomorphic to $\RR^{d}$). We use the covariant derivative point of view. We assume that the coefficients of the connections are always smooth.
\par
It is usual, and convenient, when dealing with tensors, to denote the $d$ coordinates of points and vectors in a dimension $d$ chart using exponents as in $x^1$, $x^2$, \ldots, $x^d$.
However, we chose to use $x_1$, $x_2$, \ldots, $x_d$ instead.
The reason is that, starting from Section~\ref{sec:LocalModels}, we make abundant use of quadratic forms and we prefer to write ${x_1}^2-{x_2}^2$ instead of $(x^1)^2-(x^2)^2$, and want to maintain a coherent notation throughout the article.
For all other types of tensors: forms, curvature, etc., we use the traditional notation and conventions.

\subsection{Affine connection}

Denote by $(e_{1},\dots,e_{d})$ the canonical basis of $\RR^{d}$.
Let $U$ be an open subset of $\RR^d$ considered as the domain of a chart of a manifold, and let $X$, $Y$ be vector fields on $U$.
We recall that $\partial_X Y$ denotes the \emph{directional derivative}\footnote{Recall that the vector field $\partial_X Y$ depends on the chart and will not be transported in a compatible way by a change of variables $\phi$, i.e.\ most of the time $\phi_*(\partial_X Y) \neq \partial_{\phi_*X}\phi_*Y$.} of $Y$ along $X$ in the chart $U$:
\[\partial_X Y = \sum_{i} X_{i} \frac{\partial Y}{\partial x_{i}} = \sum_{ij} X_{i} \frac{\partial Y_{j}}{\partial x_{i}} e_j.\]
In the chart, an affine connection $\nabla$ is expressed as follows: it associates to vector fields $X$ and $Y$ the vector field
$$\nabla_{X} Y = \partial_{X} Y + \Gamma(X,Y)$$
where $\Gamma$ is a bilinear map of $\RR^{d}=T_{x} U$ for each $x \in U$, taking values in $T_{x} U$:
$$\Gamma(X,Y) = \sum_{ijk} \Gamma^{i}_{jk} X_{j} Y_{k} e_{i}$$
where the coefficients $\Gamma^{i}_{jk}$ are $d^{3}$ functions of $x$ and are called the \emph{Christoffel symbols} of the connection in the chart.
By abuse of language we also call $\Gamma$ the Christoffel symbol.
The ordering of the two lower indices of $\Gamma^{i}_{jk}$ is not uniform among authors, so the formula above specifies our choice.
\par
The \emph{canonical connection} of $\RR^{d}$ is defined by $\Gamma=0$, i.e.\ 
$$\nabla_{X} Y =  \partial_{X} Y.$$
A connection on $M$ that is canonical in one chart will generally fail to be canonical in other charts.
A connection for which there are local charts where it is canonical is called \emph{locally trivializable}.

\subsubsection{Parallel transport}

For a curve $t \mapsto \gamma(t)$ whose derivative does not vanish, the parallel transport of a vector $v$ along $\gamma$ with respect to the connection $\nabla$ is a vector $v(t)$ attached to $\gamma(t)$ and which, locally, for any pair of vector fields $X$ and $Y$ satisfying $X(\gamma(t)) = \dot{\gamma}(t)$ and $Y(\gamma(t)) = v(t)$, is a solution of
$$\nabla_{X} Y = 0.$$
This amounts to the ordinary differential equation
\[\dot{v}(t) = -\Gamma\{\gamma(t)\} \big(\dot{\gamma}(t),v(t)\big)\]
which, omitting $t$, reads as
\[\dot{v} = -\Gamma\{\gamma\} \big(\dot{\gamma},v\big).\]
This formula allows us to generalize parallel transport to paths whose derivative vanishes for some values of $t$. 
The map sending $v(0)$ to $v(t)$ can be seen as a map from the tangent space of $M$ at the point $\gamma(0)$ to the tangent space of $M$ at the point represented by $\gamma(t)$. This map is also called \emph{parallel transport along $\gamma$}.

\subsubsection{Geodesics}

The famous geodesic equation for a path $\gamma$ asks that $\dot{\gamma}(t)$ be a parallel vector field along $\gamma$, i.e.
\[\ddot{\gamma} = -\Gamma\{\gamma\} \big(\dot{\gamma},\dot{\gamma}).\]
It depends only on the quadratic vector form associated with the bilinear map $\Gamma$. It has the same geodesics as the symmetrized connection with $\Gamma^{\mathrm{sym}}(u,v)=\frac{1}{2}\left(\Gamma(u,v)+\Gamma(v,u)\right)$.

\subsection{Symmetric and Levi-Civita connections}

The connection is symmetric (a.k.a.\ torsion-free) whenever the bilinear form $\Gamma$ is symmetric, i.e.\
\[\forall i,j,k,\ \Gamma^i_{jk} = \Gamma^i_{kj}.\]
This condition is independent of the chart and can be defined in a coordinate-independent way by the cancellation of an associated tensor called the torsion, but we will not use this here.
\par
Given a Riemannian metric $g$, there is a unique symmetric connection whose parallel transport preserves $g$. It is called the \emph{Levi-Civita connection} of $g$. Its parallel transport consists of isometries between tangent spaces, with respect to $g$.

\subsection{Curvature tensor}

One of the most commonly seen definitions of the \textit{curvature tensor} $R$ associated with a connection $\nabla$ is as
\[
R(X,Y)Z=\nabla_X\nabla_YZ-\nabla_Y\nabla_XZ-\nabla_{[X,Y]}Z,
\]
for all vector fields $X,Y,Z$. 
It is a $(1,3)$-tensor.
We show this definition to fix conventions but will not use it directly.
We will use that the coefficients of $R$ in charts are given (see Proposition~7.6, page~145 of Chapter~III of \cite{KN} for a reference\footnote{They decompose $R(X,Y)Z$ as $\sum R^{i}_{jkl} X_k Y_l Z_j e_i$. Hence what we denote $i,j,k,l$ they denote respectively $k,l,j,i$. Moreover in the coefficients of $R$ they put first the lower index that we put last.}) by
\begin{equation}\label{eq:CurvatureTensor}
R_{ijk}^l = \partial_i \Gamma_{jk}^{l} - \partial_j \Gamma_{ik}^{l}
+ \sum_p \left(\Gamma_{ip}^l\Gamma_{jk}^p-\Gamma_{jp}^l\Gamma_{ik}^p\right),
\end{equation}
and that they give the (opposite of the) first-order infinitesimal holonomy of the connection, as explained in Section~\ref{subsub:InfinitesimalHolonomy}.
\par
In dimension $d=1$ the curvature tensor is automatically $0$.

\begin{rem}
There are two conventions for the definition of the curvature tensor, depending on the author(s), and they give opposite tensors.
We chose the one that seems the most frequently used:
among the articles and books we cite, only \cite{Be} uses the opposite one. 
To our knowledge, all the references agree on the Ricci tensor and on the scalar curvature deduced from the latter.
\end{rem}

\subsection{Holonomy}

The parallel transport along a closed curve is called a \emph{holonomy} and is a self-map of the tangent space at $\gamma(0)$.

\subsubsection{Flat connections}

A connection is called \emph{flat}\footnote{Some authors use the term \emph{parallelizable} because of its relation with invariant frame fields. This may be less ambiguous than \emph{flat} in a situation where the connections are allowed not to be symmetric, or in situations where several notions of flatness occur (conformal flatness, projective flatness, etc.).} if, locally, its holonomies are all equal to the identity. By this we mean that every point has a neighborhood $U$ such that the holonomy of every closed path contained in $U$ is the identity map.
This is equivalent to asking that the curvature tensor vanishes everywhere (see \cite{KN}, Chapter~II, Section~9, Theorem~9.1).
\par
A connection is locally trivializable if and only if it is symmetric and flat (see \cite{Sz}, Section~18.7, page~532).
\par
In dimension $d=1$ a connection is automatically locally trivializable.

\subsubsection{Infinitesimal holonomy}\label{subsub:InfinitesimalHolonomy}

For a symmetric affine connection $\nabla$ on an open subset of $\mathbb{R}^{d}$, we consider
a closed loop $\gamma$ encompassing a small planar parallelogram directed by vectors $e_{i},e_{j}$ and positively oriented with respect to the pair $(e_i,e_j)$ of linearly independent vectors.
\par
The holonomy of the parallel transport around $\gamma$ is an endomorphism of the tangent space at $\gamma(0)$ and, as the diameter of the rectangle tends to $0$, the holonomy expands as $\Id + \mathcal{A}\cdot\mathcal{M}_{i,j} + o(\mathcal{A})$, where $\mathcal{A}$ is the area of the rectangular flag for the canonical Euclidean metric of $\RR^d$.
\par
The endomorphism $\mathcal{M}_{i,j}$ will be referred to here as the \textit{first-order infinitesimal holonomy} of the connection at this point for the rectangular flag directed by vectors $e_{i}$ and $e_{j}$.
\par
We denote by $M(i,j)$ the matrix corresponding to the endomorphism $\mathcal{M}_{i,j}$ in the canonical basis of $\RR^d$. We refer to it as a \textit{first-order infinitesimal holonomy matrix}. We will omit the mention of ``first-order'' in the sequel. The coefficient at the $l$-th row and $k$-th column of $M(i,j)$ is equal to $-R_{ijk}^{l}$, i.e.
\[ M(i,j) = (-R_{ijk}^{l})_{kl}.\]
This is, for instance, proved in \cite{Sc}, Chapter~III, Section~4, Equation~(4.11), page~140.
\par
If we replace $e_i$ and $e_j$ by any pair of non-colinear vectors $u$ and $v$ then the formula for the holonomy becomes $\Id + \mathcal{A}\cdot\mathcal{M} + o(\mathcal{A})$ with $\cal M = \sum_{i,j} \cal M_{i,j} u_i v_j$.
One of the defining properties of hybrid connections is that their first-order infinitesimal holonomies are infinitesimal isometries of the Hessian metric, see Remark~\ref{rem:equiv}.
Since infinitesimal isometries form a vector subspace of the set of endomorphisms of the tangent space at a given point, it equivalent that $M(i,j)$ is an infinitesimal holonomy matrix for all $i,j$ in any given chart.
We will see that this property can be expressed directly in terms of the coefficients of the curvature tensor.

\subsection{Change of variable}\label{sub:chvar}

If $\phi$ is a diffeomorphism between open subsets of $\RR^d$, and $\Gamma$ the Christoffel symbol of a connection in some chart, then the image of this connection under $\phi$ has Christoffel symbol $\wt \Gamma$ that is expressed as follows:
\begin{equation}\label{eq:chvar:Gamma}
\wt \Gamma\{\phi\} \circ (L,L) = L\circ\Gamma - H,
\end{equation}
where the linear map $L\in M_{d}(\mathbb{R})$ is the differential of $\phi$, $H:\RR^d \times\RR^d \to \RR^d$ is the Hessian of $\phi$ and is seen as bilinear map, $(L,L)$ is the linear map $(u,v)\mapsto (L(u),L(v))$ and $\Gamma\{\phi\}$ means that at $x\in\R^d$ we must consider $\Gamma$ at $\phi(x)$.
\par
If the change of variables is just linear: $\phi=L$, then 
\begin{equation}\label{eq:chvarlin:Gamma}
\wt \Gamma \{L\} \circ (L,L) = L\circ\Gamma.
\end{equation}
Coming back to general changes of variable $\phi$, Equation~\ref{eq:chvar:Gamma} allows us to give a simple expression of the pull-back $\nabla$ of the canonical connection $\wt\nabla$ (for which $\wt \Gamma = 0$):
\[\Gamma = L^{-1}\circ H.\]

\section{Connections with straight geodesics in flat affine manifolds}\label{sec:FlatAffine}

Two connections whose geodesics are the same as unparameterized curves are called \emph{projectively equivalent} and said to belong to the same \emph{projective class}.
A connection is called \emph{projectively flat} if it is projectively equivalent to a flat affine connection, i.e.\ if it admits charts in which all geodesics are straight lines. Such a chart is called a \emph{projective chart}.
\par
In this article we are concerned with a manifold $(M,D)$ endowed with a flat affine connection $D$ and need to determine and study the connections in the same projective class as $D$.
For a local study, it is enough to work in an affine chart of $M$, where $D$ is trivial. Hence here we determine, after Weyl, connections on open subsets of $\R^d$ whose geodesics are straight lines.
\par
In dimension one, all connections are symmetric and it is obvious that all connections belong to a unique projective class, making the theory almost pointless. However for technical reasons related to potential and Corollary~\ref{cor:aff:restr}, this helps to determine the parameterization of geodesics in any dimension, so we will also cover this case.
\par
In this section and further, we make use of the Kronecker delta: $\delta_{ij} = 0$ if $i\neq j$ and $\delta_{ij} = 1$ if $i=j$.

\subsection{The associated form}\label{sub:assoc-form}

In this section, we investigate the form that a connection takes in a projective chart.
These results are not new. In fact, Hermann Weyl already proved a more general statement, see Proposition~17, page~251, in volume~2 of \cite{Sp}.

We give a first characterization for a connection that is not necessarily symmetric.

\begin{lem}
Consider an open subset $U$ of $\mathbb{R}^{d}$ and a general connection $\nabla$ on $U$, with symbol $\Gamma$.
The following statements are equivalent:
\begin{enumerate}
\item All geodesics are straight line intervals of $\mathbb{R}^{d}$.
\item Equations~\eqref{eq-GGG} and \eqref{eq-zer} hold:
\begin{eqnarray}
\label{eq-GGG} \forall i\neq j,\ \quad \Gamma_{ii}^i & = & \Gamma^{j}_{ij}+\Gamma^{j}_{ji}\ \text{ and }\\
\label{eq-zer} \forall i,j,\ \forall k\notin\{i,j\},\ 0 & = & \Gamma_{ij}^k+\Gamma_{ji}^k\ .
\end{eqnarray}
\end{enumerate}
\end{lem}
By setting $i=j$, Equation~\eqref{eq-zer} above reads $\Gamma_{ii}^k=0$ when $k\neq i$.
Of course, when $d=1$, all the conditions of (2) are empty.

\begin{proof}
We recall that the geodesic equation reads $\gamma''(t) = -\Gamma(\gamma'(t),\gamma'(t))$. So the condition on geodesics is equivalent to
\[\forall v,\ \Gamma(v,v)\text{ is colinear to }v.\]
This reduces the question to a simple linear algebra exercise, which we now solve.

By applying this condition to $v=e_i$, we get that $\Gamma_{ii}^k=0$ for all $k\neq i$.
By applying this to $v=e_i+\lambda e_j$ when $i\neq j$, we get that
\[\Gamma_{ii} + \lambda(\Gamma_{ij}+\Gamma_{ji}) + \lambda^2 \Gamma_{jj}\]
must be colinear to $e_i+\lambda e_j$ (exponents on $\lambda$ mean powers, not indices), hence 
\[\Gamma_{ii}^k + \lambda(\Gamma_{ij}^k+\Gamma_{ji}^k) + \lambda^2 \Gamma_{jj}^k = 0\]
for all $\lambda \in \mathbb{R}$, $k \notin \lbrace{i,j\rbrace}$ and
\[\Gamma_{ii}^j + \lambda(\Gamma_{ij}^j+\Gamma_{ji}^j) + \lambda^2 \Gamma_{jj}^j = \lambda \left(\Gamma_{ii}^i + \lambda(\Gamma_{ij}^i+\Gamma_{ji}^i) + \lambda^{2} \Gamma_{jj}^i \right)\]
for all $\lambda \in \mathbb{R}$.
Then, comparing the $\lambda^{1}$ coefficients of these identities we get identities~\eqref{eq-GGG} and~\eqref{eq-zer}.

For the converse, by substituting and grouping, one eventually finds that if Equations \eqref{eq-GGG} and \eqref{eq-zer} are satisfied, we have
\[\forall v = \sum x_{i}e_i,\ \Gamma(v,v) = \left(\sum \Gamma_{ii}^i x_{i}\right) v,\]
which is colinear to $v$.
\end{proof}

For symmetric connections, these constraints are strong enough to give an explicit and simple characterization. In particular, we give a general expression for the curvature tensor $R$ associated with $\nabla$.

\begin{cor}\label{cor:expression}
Consider an open subset $U \subset \mathbb{R}^{d}$ and a \emph{symmetric} connection $\nabla$ of symbol $\Gamma$. The following statements are equivalent:
\begin{enumerate}
    \item All geodesics are straight lines of $\mathbb{R}^{d}$.
    \item There exists a $1$-form $a=\sum_{i} a_{i}dx_{i}$ such that 
\begin{equation}\label{eq:GST}
\Gamma(u,v)=a(v)u+a(u)v.
\end{equation}
\end{enumerate}
In this case, for any $i,j,k$, the coefficients of $\Gamma$ satisfy
$$
\Gamma_{ij}^{k}=\delta_{ik}a_{j}+\delta_{jk}a_{i}
$$
and $a$ is unique.
Moreover, the coefficients of the curvature tensor $R$ associated with $\nabla$ are given by
\begin{equation}\label{eq:R}
R_{ijk}^{l}
=
\delta_{il}(a_{j}a_{k}-\partial_{j}a_{k})
-
\delta_{jl}(a_{i}a_{k}-\partial_{i}a_{k})
+
\delta_{kl}(\partial_{i}a_{j}-\partial_{j}a_{i})
.
\end{equation}
\end{cor}

\begin{proof}
Adding the symmetry condition to \eqref{eq-zer} we obtain that $\Gamma_{ij}^{k}=0$ for any choice of $i,j,k$ whenever $k \notin \lbrace{ i,j \rbrace}$.
We also obtain that $\Gamma_{ii}^{i}=2\Gamma_{ij}^{j}$ for any choice of $i,j$ such that $i \neq j$. We set $a_{i}=\frac{1}{2}\Gamma_{ii}^{i}$. 
The general formula $\Gamma_{ij}^{k} = \delta_{ik}a_{j} + \delta_{jk}a_{i}$ and the expression of $\Gamma(u,v)$ given in the statement immediately follow.
\par
Conversely, a connection as in (2) satisfies the relations of~\eqref{eq-zer} and since $\Gamma_{ii}^i=2a_i$, the coefficients $a_i$ are unique.
\par
The coefficients of $R$ in coordinates are given by \eqref{eq:CurvatureTensor}. We obtain the following expression for $\Gamma_{ip}^{l}\Gamma_{jk}^{p} - \Gamma_{jp}^{l}\Gamma_{ik}^{p}$:
$$
(\delta_{il}\delta_{jp}-\delta_{jl}\delta_{ip})a_{p}a_{k}
+(\delta_{jp}a_{i}-\delta_{ip}a_{j})\delta_{pl}a_{k}
+(\delta_{il}a_{j}-\delta_{jl}a_{i})\delta_{kp}a_{p}.
$$
Summing over every index $p$, the last two terms cancel and we obtain:
$$
\sum_{p} 
\Gamma_{ip}^{l}\Gamma_{jk}^{p} - \Gamma_{jp}^{l}\Gamma_{ik}^{p}
= \delta_{il}a_{j}a_{k} - \delta_{jl}a_{i}a_{k}
.
$$
The final formula follows.
\end{proof}

Expressed in terms of operators, Equation~\ref{eq:R} reads:
\begin{equation}\label{eq:R2}
\begin{split}
R(X,Y)Z & = a(Z)(a(Y)X-a(X)Y) 
\\
& + (\partial_Y a)(Z) X - (\partial_X a)(Z) Y
\\
& - da(X,Y) Z
\end{split}
\end{equation}
where $\partial_X a $ denotes the directional derivative of $a$ in a chart, i.e.\ $(\partial_X a) (Y) = \sum_{i,j} \partial_i a X_i Y_j$, and $da$ is the exterior derivative of $a$, which can also be expressed as $da(X,Y)=(\partial_Y a)(X)-(\partial_X a)(Y)$.

In dimension $d=1$, the conditions above are always satisfied and the only coefficient of $a$ and the only coefficient of $\Gamma$ are related by $2a_1 = \Gamma_{11}^1$. Moreover, the curvature tensor is always $0$.

\begin{defn}\label{defn:Asso}
In the present article, the differential form
\[a=\sum a_i\, dx_{i}\]
appearing in the previous statement will be called the
\emph{associated form}.
\end{defn}

In Section~\ref{sub:aff:chvar}, we prove that $a$ is independent,
as a differential form, of the affine chart of $M$ in which we work.
Note that a projective change of variables $\phi$ does not change the projective class but changes the differential form $a$ if $\phi$ is not affine.

\subsection{Ricci tensor and incompressible connections in affine charts}

The Ricci tensor $\Ric$ is a $(0,2)$-tensor obtained by contraction of the curvature tensor $R$. Its coefficients are given (with our choices of conventions for $R$) by 
\begin{equation}\label{eq:defRic}
\Ric_{jk}=\sum_{i} R_{ijk}^{i}.
\end{equation}

Thus, in an affine chart of a $d$-dimensional affine manifold $(M,D)$, the coefficients of the Ricci tensor of a symmetric affine connection $\nabla$ deduced from Equation~\eqref{eq:R} are written as follows, in terms of 
the coefficients of the associated form $a$ and their partial derivatives:
\begin{equation}\label{eq:Ric}
\Ric_{jk}=(d-1)a_{j}a_{k}+\partial_{k}a_{j}-d\,\partial_{j}a_{k}.
\end{equation}
Note that in this equation, $d$ refers to the dimension of the manifold, and the last term is the product of the real numbers $d$ and $\partial_{j}a_{k}$, not the differential of $\partial_{j}a_{k}$.

Below we prove that a projectively flat connection is flat if and only if it is Ricci flat.

\begin{prop}\label{prop:flatequiv}
Given a $d$-dimensional affine manifold $(M,D)$ and a symmetric affine connection $\nabla$ that is projectively equivalent to $D$, the following statements are equivalent:
\begin{enumerate}
    \item $\nabla$ is Ricci flat (i.e.\ $\Ric=0$);
    \item $\nabla$ is flat (i.e.\ $R=0$).
\end{enumerate}
\end{prop}

\begin{proof}
If $d=1$, then $R=0$ and $\Ric=0$ so it is enough to handle the case $d\geq 2$.
Since $R$ and $\Ric$ are tensors, if they vanish in a chart then they vanish for any other chart $\tilde \phi$ on the corresponding points.
\newline
If $R$ vanishes, then $\Ric$ vanishes according to Equation~\eqref{eq:defRic}.
\newline
Conversely, assume that $\Ric$ vanishes. In particular, $\Ric$ is symmetric, which, by Equation~\eqref{eq:Ric}, means that for any pair of indices $i,j$,
$(d-1)a_{i}a_{j}+\partial_{j}a_{i}-d\,\partial_{i}a_{j} = 
(d-1)a_{j}a_{i}+\partial_{i}a_{j}-d\,\partial_{j}a_{i}$ which, after simplifications, leads to $\partial_{i}a_{j}=\partial_{j}a_{i}$. Substituting this into Equation~\eqref{eq:Ric}, the condition $\Ric_{ij}=0$ then simplifies to $a_{i}a_{j}=\partial_{i}a_{j}$.
Injecting these identities into Equation~\eqref{eq:R}, we obtain that every coefficient of the curvature tensor $R$ vanishes.
\end{proof}

The incompressibility property (see Definition~\ref{defn:Incompressible}) is equivalent to the infinitesimal holonomy matrices (see Section~\ref{subsub:InfinitesimalHolonomy}) having vanishing trace everywhere, i.e.\ $\sum_{k} R_{ijk}^{k} = 0$.

\begin{prop}\label{prop:Ricci}
Given a $d$-dimensional affine manifold $(M,D)$, we consider an affine chart $\phi$.
Given a symmetric connection $\nabla$ that is projectively equivalent to $D$, let $a_{1},\dots,a_{d}$ be the coefficients of the form $a$ associated with $\nabla$ in chart $\phi$. The following statements are equivalent:
\begin{enumerate}
    \item The associated form $a$ is closed.
    \item There exists a locally defined nowhere vanishing smooth function $\psi$, called a \emph{logarithmic potential}, such that for any index $i$, $a_{i}=-\frac{\partial_{i}\psi}{\psi}$ (i.e.\ $a=-\frac{d \psi}{\psi}$).
    \item The Ricci tensor of $\nabla$ is symmetric ($\Ric_{jk}=\Ric_{kj}$ for any pair of indices $j,k$).
    \item The connection $\nabla$ is incompressible.
\end{enumerate}
In these cases, the logarithmic potential is unique up to multiplication by an element of $\RR^*$, the preserved volume forms are locally given by
\[\omega = \alpha \frac{dx_{1}\wedge\cdots\wedge dx_{d}}{\psi^{d+1}}\]
for $\alpha\in\RR^{\ast}$,
and the coefficients of the Ricci tensor are given by
\begin{equation}\label{eq:Ric-psi}
\Ric_{jk}=(d-1)\frac{\partial_{j}\partial_{k}\psi}{\psi}.
\end{equation}
Moreover, the coefficients of the curvature tensor are completely determined by the Ricci tensor: if $d=1$ then $R=0$ and $\Ric =0$; if $d\geq 2$ then
$$
R_{ijk}^{l}=\frac{1}{d-1}(\delta_{il}\Ric_{jk}-\delta_{jl}\Ric_{ik})
.$$
We can also express $R$ in terms of $\psi$: 
\begin{equation}\label{eq:Rpsi}
R_{ijk}^{l} = \delta_{il}\frac{\partial_j\partial_k\psi}{\psi} - \delta_{jl}\frac{\partial_i\partial_k\psi}{\psi}.
\end{equation}
\end{prop}

\begin{proof}
By Equation~\eqref{eq:Ric} we have
\[
\Ric_{jk}-\Ric_{kj} = (d+1) (\partial_k a_j - \partial_j a_k)
.\]
This proves that the Ricci tensor is symmetric if and only if the differential form $a=\sum_{i} a_{i}dx_{i}$ is closed ($\partial_{k}a_{j}=\partial_{j}a_{k}$ is satisfied for any pair of indices $j,k$).
\par
The closedness of $a$ is equivalent to the existence of a locally defined function $\theta$ such that $a=d\theta$. We define $\psi=e^{-\theta}$ and obtain that $a=d\theta=-\frac{d\psi}{\psi}$. 
Conversely, if (2) holds then $a=d\theta$ for $\theta = -\log |\psi|$.
We proved that statements (1), (2) and (3) are equivalent.
\par
The connection $\nabla$ locally preserves some volume form if and only if its first-order infinitesimal holonomy matrices are traceless (see Section~\ref{subsub:InfinitesimalHolonomy}).
This is equivalent to saying that for any pair of indices $i,j$, we have $\sum_{k} R_{ijk}^{k}=0$. By formula~\eqref{eq:R}, we have
$$
\sum_{k} R_{ijk}^{k} = (d+1)(\partial_{i}a_{j} - \partial_{j}a_{i})
.
$$
Statement (4) is thus equivalent to the closedness of differential form $a$.

Near a given point, an invariant volume form is uniquely determined by its value at that point. Indeed, any prescribed value extends uniquely to a volume form on a sufficiently small neighborhood by parallel transport, and incompressibility guarantees that this extension is independent of the chosen path. Since the space of volume forms at a point is one-dimensional, the correspondence between the initial value and the resulting local invariant volume form is linear, and every local invariant volume form is obtained in this way.
\par
The analysis of the parallel transport will naturally lead us to the formula $\frac{dx_{1}\wedge\cdots\wedge dx_{d}}{\psi^{d+1}}$.
Indeed, consider an initial point $x_*\in \RR^d$, choose one of the basis vectors $e_i$, and define the path $\gamma(t) = x_*+te_i$.
Let $v(0)$ be any vector and $v(t)$ be the parallel transport of $v$ along $\gamma$.
By definition $\dot v = -\Gamma(\dot\gamma,v) = \frac{1}{\psi} \left(d\psi(\dot\gamma) v + d\psi(v) \dot\gamma\right)$
where $\Gamma$, $\psi$ and $d\psi$ are evaluated at the point $\gamma(t)$ and for any function $f$ of $t$, $\dot f$ means $f'(t)$.
We obtain
\[\dot v = \frac{1}{\psi} \left(d\psi(e_i) v + d\psi(v) e_i\right).\]
Note that $d\psi(e_i) = \dot\psi = (\psi\circ\gamma)'(t)$.

We choose $v(0)=e_j$. Then, $v(t)$ remains in the vector subspace generated by $e_i$ and $e_j$.
\par
If $j\neq i$ then let us write $v(t) = f(t)e_i+g(t) e_j$ with $f(0)=0$ and $g(0)=1$. Then the differential equation reads
$\dot f e_i + \dot g e_j = (\cdots)e_i + \frac{\dot\psi}{\psi}g e_j$. It follows that $\dot g/g = \dot\psi /\psi$ where defined, hence $g(t) = \psi\circ\gamma(t)/\psi(x_*)$.
\par
If $j=i$ then let us write $v(t) = f(t) e_i$ with $f(0)=1$ and the differential equation reads
$\dot f e_i = 2\frac{d\psi(e_i)}{\psi} f e_i = \frac{2\dot\psi}{\psi} f e_i$ where defined.
Hence $\dot f/f = 2\dot \psi/\psi$, so $f(t) = (\psi\circ\gamma(t)/\psi(x_*))^2$.
\par
It follows that the form $\frac{dx_{1}\wedge\cdots\wedge dx_{d}}{\psi^{d+1}}$ is invariant under transport along lines parallel to $e_i$.
And this holds for any $i$. Any volume form which is invariant by parallel transport must be invariant by transport along these lines, so must be proportional to the above form. Since we know by incompressibility that there is an invariant form, the above form is invariant.
\par
The last expressions for $\Ric$ and $R$ are obtained by simple computations, using:  Equation~\eqref{eq:Ric} and point (2) of the present proposition for $\Ric$; Equations~\eqref{eq:R} and~\eqref{eq:Ric} and the fact that $a$ is closed for $R$.
\end{proof}
Actually, (3) $\iff$ (4) holds for a general connection, by \cite{NSAff}, Chapter~I, Proposition~3.1.

\begin{rem}
Note that the potential of $\nabla$ is non-unique and a priori only defined locally.
If $U\subset\R^d$ is not simply connected, it may happen that $\psi$ cannot be globally defined on $U$, in the sense that there would be monodromy factors to take into account.
For instance, let $U=\R^2 - \lbrace{ O \rbrace}$, where $O$ is the origin, and $\psi = \exp(\theta)$ where $\theta$ is the polar angle of vector $(x,y)$.
Here the monodromy factor is $\exp(2\pi)$.
\end{rem}

Potentials behave like functions under a change of affine charts, see Section~\ref{sub:aff:chvar}.

\subsection{Restrictions and affine changes of variables}\label{sub:aff:chvar}

Here we explore the effect of restrictions to affine subspaces, and of affine changes of variables, in affine charts of $(M,D)$, for a symmetric connection projectively equivalent to $D$, on the potential and the associated form.

A remarkable property of the connections appearing in Corollary~\ref{cor:expression} is that parallel transport preserves planes.

\begin{cor}\label{cor:plane}
Assume that $d\geq 2$ and $U\subset \RR^d$ is a chart in which a symmetric connection $\nabla$ has its geodesics that are straight lines.
Consider a path $\gamma$ of $U$ that is also contained in an affine $2$-plane $P$ of $\RR^d$. If $v$ is a vector contained in the vector subspace $\vv P$ of $\RR^d$ associated with $P$, then the parallel transport of $v$ along $\gamma$ still belongs to $\vv P$.
\end{cor}

\begin{proof}
From the formula $\Gamma(u,v) = a(u)v+a(v)u$ it follows that $\Gamma(u,v)$ belongs to the vector space generated by $u$ and $v$.
Recall the parallel transport equation: $\dot v= -\Gamma (\dot\gamma,v)$ where $\Gamma$ is taken at point $\gamma(t)$.
It follows that if $v$ and $\dot\gamma$ belong to $\vv P$ then $\dot v$ does too.
The restriction of the connection $\nabla$ to the affine plane $P$ of $\RR^d$ thus makes sense, and defines a parallel transport $v(t)$ of $v$ along $\gamma$, for which $v(t)\in\vv P$, and which satisfies the same equation $\dot v= -\Gamma (\dot\gamma,v)$.
The result follows by uniqueness of solutions of differential equations.
\end{proof}

In an affine manifold, the Christoffel symbol $\Gamma$ is a well-defined (1,2)-tensor, i.e.\ its expression in various affine charts differs by the same formula as for a (1,2)-tensor (see Equation~\eqref{eq:chvarlin:Gamma} in Section~\ref{sub:chvar}).
At a point $m\in M$, it gives a well-defined bilinear map $\Gamma\{m\} : T_m M\times T_m M\to T_m M$.

Consider an affine manifold $(M,D)$.
An affine submanifold $N$ is a submanifold such that $N$ is, in any chart, locally an affine subspace of $\RR^d$.
We saw in Corollary~\ref{cor:plane} that for a connection in $\R^d$ that is projectively equivalent to the canonical connection, the parallel transport along an affine $2$-plane $P$ of $\R^d$ preserves vectors tangent to $P$.
This generalizes to any affine submanifold of $M$.

\begin{cor}\label{cor:aff:restr}
Given an affine manifold $(M,D)$, we consider a symmetric affine connection $\nabla$ that is projectively equivalent to $D$ and denote $\Gamma$ its Christoffel symbol in affine charts. For any affine submanifold $N$ of $M$:
\begin{itemize}
\item for any point $m\in M$, $\Gamma\{m\}(T_m N,T_m N) \subset T_m N$;
\item $\nabla$ has a well-defined restriction $\nabla^N$ to $N$, given by the restriction of $\Gamma$;
\item the restriction of the form $a$ to $N$ is the associated form of $\nabla^N$;
\item the parallel transport of any vector in $T N$ along any curve contained in $N$ remains in $T N$;
\item $\nabla^N$ is projectively equivalent to $D$.
\end{itemize}
Assume moreover that $\nabla$ is incompressible. Then
\begin{itemize}
\item $\nabla^N$ is incompressible;
\item for any logarithmic potential $\psi$ for $\nabla$ (see Proposition~\ref{prop:Ricci}), defined on an open subset $V$ of $M$, the restriction of $\psi$ to $N\cap V$ is a logarithmic potential for $\nabla^N$ on $N\cap V$.
\end{itemize}
\end{cor}
\begin{proof}
Let us work in a chart.
By Equation~\ref{eq:GST}, $\Gamma(u,v) = a(u)v+a(v)u$.
The first three points immediately follow.
The fourth point is a consequence of uniqueness of solutions of differential equations.
The connections $\nabla$ and $\nabla^N$ having the same geodesics in $N$, the fifth point is immediate.
If $\nabla$ is incompressible then $a$ is closed, so the restriction of $a$ to $N$ is closed, so it is incompressible (see Proposition~\ref{prop:Ricci}).
The relation $a=-d\psi/\psi$ still holds for the restrictions of $a$ and $\psi$ to $N$, so by Proposition~\ref{prop:Ricci}, the last point follows.
\end{proof}

\begin{prop}
Consider an affine manifold $(M,D)$ and an affine chart $\phi : O\subset M \to U\subset\RR^d$.
Consider a symmetric connection $\nabla$ on $M$
with associated form $a = \sum_{i} a_{i}dx_{i}$ in chart $\phi$.
If $A$ is an affine change of variables then $A\circ\phi$ is also an affine chart.
Denote by $L$ the linear part of $A$ and $\tilde{a} = \sum_{i}\tilde{a}_{i}dx_{i}$ the associated form for this new chart.
We have:
\begin{equation}\label{eq:chvar:real:aff}
\tilde{a}\lbrace{ A \rbrace}= a \circ L^{-1}
.
\end{equation}
This notation is a shorthand for the following equality: for every $x\in U$ and every $u \in \mathbb{R}^{d}$,
$\tilde{a}\{A(x)\}(u) = a\{x\}(L^{-1}(u))$ where $a\{x\}(u) = \sum_{i} a_{i}(x)u_{i}$ and 
$\tilde{a}\{A(x)\}(u) = \sum_{i} \tilde{a}_{i}(A(x))u_{i}$.
\par
Moreover, if $\nabla$ is incompressible, then the potential $\tilde \psi$ in the new chart and the potential $\psi$ in the original chart are related by:
\[ \psi = \tilde\psi \circ A
.\]
\end{prop}

\begin{proof}
This is an easy consequence of Equations~\eqref{eq:chvarlin:Gamma} and \eqref{eq:GST}.
Indeed, consider the first one: $\tilde \Gamma \{A\} \circ (L,L) = L\circ \Gamma$.
Substitute the second equation on both sides:
\[
(\tilde a\{A\}\circ L(v)) L(u) + (\tilde{a}\{A\}(u)) L(v) =
L(a(v)u+a(u)v) = a(v)L(u)+a(u)L(v).
\]
By composing on the right with $(L^{-1},L^{-1})$, we obtain that 
\[
(\tilde a\{A\}(v)) u + (\tilde{a}\{A\}(u)) v = (a\circ L^{-1}(v))u+(a\circ L^{-1}(u))v.
\]
for all $u,v$. By uniqueness of the $a_i$ in Corollary~\ref{cor:expression}, it follows that $\tilde a\{A\} = a\circ L^{-1}$.
\end{proof}

\begin{rem}
Note that \eqref{eq:chvar:real:aff} is exactly the same formula as for differential forms. In other words, $a$ behaves as a $1$-form on $(M,D)$. 
\end{rem}

\section{Hessian metrics and characterization of hybrid connections}\label{sec:HessianMetrics}

\subsection{Hessian operator}

On a real manifold $M$ endowed with a flat affine connection $D$, the Hessian operator associates to any smooth function $f \in \mathcal{C}^{\infty}(M,\mathbb{R})$ a field $H(f)$ of symmetric bilinear forms whose matrices
are $\frac{\partial^{2}f}{\partial x_{i}\partial x_{j}}$ in a local system of affine coordinates $x_{1},\dots,x_{d}$.
\par
As a $(0,2)$-tensor, it is independent of the affine chart.\footnote{This follows immediately from the fact that changes of affine coordinates are linear. Note that for any connection, one can define an associated Hessian operator, but this is out of the scope of this article.}
\par
On a flat affine manifold $(M,D)$, a pseudo-Riemannian metric $g$ is \textit{Hessian} if it locally coincides with $H(f)$ for some locally defined smooth potential $f$. Such a potential $f$ is locally unique up to addition of an affine function with respect to $D$.

\begin{ex}
On $\mathbb{R}^{d}$ endowed with its usual affine structure $D$, the quadratic potential $\psi(x) = \frac{1}{2} \sum\limits_{i \leq p}{x_{i}}^{2} - \frac{1}{2} \sum\limits_{j>p} {x_{j}}^{2}$ induces the standard pseudo-Euclidean metric $g$ of signature $(p,d-p)$. The triple $(\mathbb{R}^{d},D,g)$ is a Hessian manifold.
\end{ex}

\begin{rem}\label{rem:k1}
Let us mention that neither the Poincaré disk or half-plane models of the hyperbolic plane, nor its Cayley-Klein disk model, are Hessian with respect to the canonical connection of $\R^2$.
According to \cite{bryant2025hessianizability}, every pseudo-Riemannian manifold of dimension $2$ has local charts where it is Hessian.
Actually, the whole hyperbolic plane has a Hessian model\footnote{Communicated to us by Gabriel Khan, attributed to Amari.} over the Euclidean upper half plane, of potential $\psi = \frac{x^2}{y} - \log y$.
\end{rem}

\subsection{Infinitesimally conformal holonomy and coefficients of the curvature tensor}\label{sub:InfinitesimalHolonomy}

We say that a linear map $\phi$ from a pseudo-Euclidean vector space $E$ to itself is a \textit{similarity} if there exists a \emph{positive} constant $k \in \mathbb{R}_{>0}$ such that for any $v \in E$, $q\circ \phi(v) = k\cdot q(v)$. We say that such a map is conformal if moreover its determinant is positive.
These two classes of maps define closed Lie subgroups of $\on{GL}(E)$ and their tangent space at the identity (their Lie algebra) is the same. Their elements will be referred to interchangeably as \emph{infinitesimal similarities} or \emph{infinitesimally conformal maps}.

\begin{prop}\label{prop:InfinitesimalHolonomy}
Consider a Hessian manifold $(M,D,g)$ and a symmetric connection $\nabla$ that is projectively equivalent to $D$.
For any point $m \in M$, the infinitesimal holonomy of $\nabla$ at $m$ is an infinitesimal similarity with respect to the pseudo-Riemannian metric $g_m$ ($g$ at point $m$) if and only if in the affine chart of $(M,D)$ where $g_m$ writes as $\sum\limits_{i=1}^{d} \epsilon_{i}\, d{x_{i}}^{2}$ for some $\epsilon_{1},\dots,\epsilon_{d} \in \lbrace{ 1,-1 \rbrace}$, the coefficients of its curvature tensor $R$ at point $m$ satisfy the following conditions:
\begin{enumerate}
    \item for any indices $i,j,k,l$ such that $k \neq l$, we have $R_{ijk}^{l}=-\epsilon_{l}\epsilon_{k} R_{ijl}^{k}$ (where $\epsilon_{1},\dots,\epsilon_{p}=1$ while $\epsilon_{p+1},\dots,\epsilon_{d}=-1$) at point $m$;
    \item for any indices $i,j$, $R_{ijk}^{k}$ at point $m$ does not depend on $k$.
\end{enumerate}
The first-order holonomy of $\nabla$ consists of infinitesimal isometries if and only if we have that (1) holds and that for any indices $i,j,k$, $R_{ijk}^{k}=0$ at point $m$.
\end{prop}

\begin{proof}
In the canonical basis of $\mathbb{R}^{d}$, we denote by $M(i,j)$ the holonomy matrix for an infinitesimal rectangle based at some point and directed by the vectors $e_{i}$ and $e_{j}$. The matrix $M(i,j)$ has its coefficient at the $l$-th row and $k$-th column equal to $-R_{ijk}^{l}$.\newline
In an affine chart where $g_m$ writes as $\sum_{i=1}^{d} \epsilon_{i}\, d{x_{i}}^{2}$, the matrices $M(i,j)$ are infinitesimal similarities of $g_m$ if and only if they can be written as the sum of a scalar multiple of the identity matrix and of a matrix $A$ such that $I_{p,d-p}A$ is antisymmetric ($I_{p,d-p}$ is the diagonal matrix whose $p$ first coefficients are equal to $1$ while the remaining $d-p$ are equal to $-1$). The condition on the coefficients of the curvature tensor follows.
\par
Infinitesimal isometries correspond to the fact that the diagonal terms vanish.
\end{proof}

\subsection{Characterization of hybrid connections}\label{sub:charHybrid} 

By Proposition~\ref{prop:Ricci}, the coefficients of a symmetric and incompressible connection $\nabla$ projectively equivalent to $D$ depend only on a logarithmic potential $\psi$.
In the lemma below, we translate the infinitesimal conformality condition into a system of partial differential equations satisfied by the potential.

\begin{lem}\label{lem:HessianInfinitesimal0}
Consider a symmetric and incompressible connection $\nabla$ on an open subset $U$ of $\R^d$. Assume $\nabla$ is projectively equivalent to the trivial connection of $\R^d$, i.e.\ that its geodesics are supported by straight lines of $\R^d$.
Let $m\in U$ and let $\psi$ be a logarithmic potential for $\nabla$ defined in a neighborhood of $m$.
Consider a pseudo-Euclidean metric $g_m$ at $m$, i.e.\ a non-degenerate quadratic form on $T_m\R^d$. 
Then the first-order infinitesimal holonomies of $\nabla$ at $m$ are infinitesimal isometries of $g_m$ if and only if there exists a real number $\lambda$ such that $H(\psi)$ at $m$ is equal to $\lambda \cdot g_m$.
\end{lem}

\begin{proof}
Up to an affine change of variables in $\R^d$, we can put $g$ into the canonical form
$$g_m = \sum\limits_{i=1}^{d} \epsilon_{i}\, d{x_{i}}^{2}$$
for some $\epsilon_{1},\dots,\epsilon_{d} \in \lbrace{ 1,-1 \rbrace}$ satisfying 
$\epsilon_{1},\dots,\epsilon_{p}=1$ while $\epsilon_{p+1},\dots,\epsilon_{d}=-1$.
\par
Assume the first-order infinitesimal holonomies of $\nabla$ at $m$ are infinitesimal isometries of $g_m$.
By the direct implication in Proposition~\ref{prop:InfinitesimalHolonomy}:
\begin{enumerate}
    \item for any indices $i,j,k,l$ such that $k \neq l$, we have $\epsilon_{k} R_{ijk}^{l}+\epsilon_{l} R_{ijl}^{k}=0$ at $m$
    \item for any indices $i,j,k$, $R_{ijk}^{k}=0$ at $m$.
\end{enumerate}
Let us recall Equation~\eqref{eq:Rpsi} in Proposition~\ref{prop:Ricci}:
\[R_{ijk}^{l} = \delta_{il}\frac{\partial_j\partial_k\psi}{\psi} - \delta_{jl}\frac{\partial_i\partial_k\psi}{\psi}.\]
We deduce from the second point above that for any indices $i,j,k$, we have $\delta_{jk} \partial_{i}\partial_{k}\psi= \delta_{ik} \partial_{j}\partial_{k}\psi$ at $m$.
Therefore, for any $i \neq j$, we have $\partial_{i}\partial_{j}\psi=0$ at $m$ (by taking $k=j$).
\par
Now apply the first equation to $k=j$ and $l=i$ for $i\neq j$.
We get $\epsilon_j R_{ijj}^{i}+\epsilon_{i} R_{iji}^{j}=0$ at $m$.
Using the expression of $R_{ijk}^l$ given above, we get $\epsilon_{i}\partial_i^2\psi=\epsilon_{j}\partial_j^2\psi$ at $m$ for any pair of distinct indices $i,j$. In other words, there is a real number $\lambda$ such that $H(\psi)=\lambda \cdot g$ at $m$.
\par
Conversely, if $H(\psi)=\lambda g_m$ for some $\lambda\in\R$ then for all $i,j$, $\partial_i\partial_j \psi=\lambda \delta_{ij}\eps_j$ at $m$.
We get by Equation~\eqref{eq:Rpsi} that $R_{ijk}^l = \frac{\lambda}{\psi}(\delta_{il}\delta_{jk}-\delta_{jl}\delta_{ik})\eps_k$
so $R_{ijk}^k = \frac{\lambda}{\psi}(\delta_{ik}\delta_{jk}-\delta_{jk}\delta_{ik})\eps_k=0$, so~(2) holds, and $\epsilon_{k} R_{ijk}^{l}+\epsilon_{l} R_{ijl}^{k} = \frac{\lambda}{\psi}(\delta_{il}\delta_{jk} - \delta_{jl}\delta_{ik} + \delta_{ik}\delta_{jl} - \delta_{jk}\delta_{il}) = 0$ so (1) holds, hence by the converse implication in Proposition~\ref{prop:InfinitesimalHolonomy}, the first-order infinitesimal holonomy of $\nabla$ at $m$ preserves $g_m$.
\end{proof}

In order to prove Theorem~\ref{thm:MAIN}, it remains essentially to prove that the coefficient $\lambda$ of Lemma~\ref{lem:HessianInfinitesimal0} is constant, and to check we can get a globally defined hybrid potential $\psi$ of $\nabla$.

\begin{proof}[Proof of Theorem~\ref{thm:MAIN}]
We first argue locally. Let $m\in M$ and consider a logarithmic potential $\psi$ for $\nabla$ defined on a connected open neighborhood $U$ of $m$.
By Lemma~\ref{lem:HessianInfinitesimal0}, $H(\psi)=\lambda(x)\cdot g$ on the domain of $\psi$, for some function $x\in U\mapsto \lambda(x)\in\R$.
Since $g$ depends smoothly on $x$ and is non-degenerate, we deduce that $\lambda$ is a smooth function.
By hypothesis, the Riemann curvature $R$ of $\nabla$ does not vanish; hence, neither does the Ricci tensor $\Ric$ by Proposition~\ref{prop:flatequiv}.
By Equation~\eqref{eq:Ric-psi} in Proposition~\ref{prop:Ricci}, $\Ric_{jk}=(d-1)\frac{\partial_{j}\partial_{k}\psi}{\psi}$.
Hence $H(\psi)\neq 0$.
We deduce that $\lambda(x)\neq 0$.
\par
We prove now that $d\lambda=0$ at $m$.
Consider a local Hessian potential $f$ for $g$.
There is an affine chart (with respect to the flat affine connection $D$) where the matrix of $g$, which is equal to the hessian matrix $(\partial_{i}\partial_{j} f)_{i,j}$ of $f$, is diagonal at the point $m$.
Therefore the Hessian matrix of $\psi$ is diagonal too in this chart at $m$.
The identity $H(\psi)=\lambda \cdot g$ reads
$\partial_{j}\partial_{k}\psi=\lambda \partial_{j}\partial_{k}f$ in a neighborhood of $m$.
Differentiating, we get that for any indices $i,j,k$, $\partial_{i}\partial_{j}\partial_{k}\psi=(\partial_{i}\lambda)\partial_{j}\partial_{k}f + \lambda \partial_{i}\partial_{j}\partial_{k} f$. Therefore, the product $(\partial_{i}\lambda)\partial_{j}\partial_{k}f$ is invariant under permutation of indices.
For $k=j$, we deduce that $(\partial_{i}\lambda)\partial_{j}\partial_{j}f = (\partial_{j}\lambda)\partial_{i}\partial_{j}f$.
For $i \neq j$, the right hand side vanishes at $m$ by our choice of affine chart.
Moreover, nondegeneracy of $H(f)$ implies that $\partial_{j}\partial_{j}f$ does not vanish at $m$. It follows that $\partial_{i}\lambda$ vanishes at $m$ for any index $i$. Therefore, $d\lambda$ vanishes at $m$. Since this holds for every point $m\in U$ we get that $\lambda$ is locally constant and since $U$ is connected, $\lambda$ is constant.
\par
We have shown that $H(\psi)=\lambda\cdot g$ for a constant $\lambda\in\R^*$.
It follows that, $H(\psi_0) = g$ for $\psi_0 = \psi/\lambda$.
The logarithmic potential is locally uniquely defined up to a multiplication by a non-zero constant, so $\psi_0$ is also a potential for $\nabla$.
Among the potentials for $\nabla$ on $U$, $\psi_0$ is the unique one such that $H(\psi_0)=g$.
\par
We proved that for any point $m \in M$, there is a neighborhood $U$ of $m$ on which $\nabla$ has a unique\footnote{Another way to see uniqueness is that if there were two distinct potentials $\psi_{1}$ and $\psi_{2}$ that were at the same time a logarithmic potential for $\nabla$ and a Hessian potential for $g$, then their quotient $\frac{\psi_{1}}{\psi_{2}}$ would be equal to a constant $A$ while their difference $\psi_{1}-\psi_{2}$ would be an affine function $\sigma$ (with respect to the affine structure defined by $D$).
It would follow that $\psi_{1}-\psi_{2}=(1-A)\psi_{2}=\sigma$ and that $\psi_{2}$ is an affine function, so $g=H(\psi_{2})=0$. This is impossible.} logarithmic potential $\psi_0$  such that $H(\psi_0)=g$ on $U$.
It follows that $\psi_0$ is globally defined.
\par
Conversely, if the Hessian metric $g$ has a globally defined potential $\psi$ that is also a logarithmic potential for $\nabla$ in the affine charts of $M$, then by the converse in Proposition~\ref{lem:HessianInfinitesimal0}, the first-order infinitesimal holonomies of $\nabla$ are infinitesimal isometries of $g$.
Moreover we can reverse the arguments of the first paragraph: since $\frac{\psi}{d-1}\Ric = H(\psi) =g$ does not vanish, Proposition~\ref{prop:flatequiv} ensures that the curvature tensor of $\nabla$ does not vanish.
Being defined by a potential, $\nabla$ is incompressible.
We have checked that $\nabla$ is hybrid.
\end{proof}

As an immediate consequence of the formula expressing $\Ric_{ij} = (d-1)\frac{\partial_i \partial_j \psi}{\psi}$ in Proposition~\ref{prop:Ricci}:

\begin{cor}\label{cor:HybridRic}
Given a hybrid connection $\nabla$ defined by a hybrid potential $\psi$ on a Hessian manifold $(M,D,g)$, the Ricci tensor $\Ric$ associated with $\nabla$ is in the conformal class of the pseudo-Riemannian metric $g$ where $\psi>0$ and of $-g$ where $\psi<0$:
$$
\Ric=\frac{d-1}{\psi}g.
$$
\end{cor}

Under the conditions of Lemma~\ref{lem:HessianInfinitesimal0}, $H(\psi)$ at $m$ is either $0$ if $\lambda=0$ or nondegenerate if $\lambda\neq 0$. Hence another consequence of this lemma is:

\begin{cor}\label{cor:slcmgk}
Consider an affine manifold $(M,D)$ and a symmetric connection $\nabla$ that is projectively equivalent to $D$. Assume that $H(\psi)$ is non-degenerate for the local logarithmic potentials $\psi$ of $\nabla$.
Then the pseudo-Riemannian metrics that are invariant under the infinitesimal holonomies of $\nabla$ are locally given by $g_m=f(m) H(\psi)$, where $f$ is a smooth non-vanishing function.
\end{cor}

In other words, it is the set of metrics that are conformal to $g=H(\psi)$ or its opposite.
On the opposite side, if $H(\psi)$ vanishes on $M$, then any metric is preserved by the infinitesimal holonomies of $\nabla$, which are all trivial since $\nabla = D$. 

\begin{rem}
Is every metric in the conformal class of the Hessian relative to some flat connection $D$ of some potential $\psi$?
In dimension $>2$, this does not hold, as a parameter counting argument shows (communicated to us by Gabriel Khan).
We already mentioned in \ref{rem:k1} that in dimension $2$, Robert Bryant proved something stronger: any metric is locally \emph{equal} to the Hessian relative to some flat connection $D$ of some potential $\psi$.
\end{rem}

\subsection{Scalar curvature of a hybrid connection}\label{sub:ScalarCurvature}

For a pseudo-Riemannian metric, one of several equivalent definitions of scalar curvature is as the trace of the $(1,1)$-tensor obtained from the $(0,2)$-tensor $\Ric$ via the musical isomorphism between $TM$ and $T^{\ast}M$.
One can raise either the first or second variable: even though the two operators thus obtained may differ, they have the same trace.
For an arbitrary (symmetric) connection, there is no notion of scalar curvature, because even though it has a Ricci tensor, there is no canonical choice of a musical isomorphism.
\par
In our particular situation of a hybrid connection, however, we can use the pseudo-Riemannian metric $g$ to define the isomorphism.
With this convention, we get what we call the \emph{scalar curvature $S(\nabla,g)$ of $\nabla$ with respect to $g$}. From Corollary~\ref{cor:HybridRic} the two endomorphisms whose trace can be taken are both equal to $\frac{d-1}{\psi} \on{Id}$ (they are equal because $\Ric$ is symmetric in our case), so
\[
S(\nabla,g) = \frac{d(d-1)}{\psi}
.\]

\begin{rem}\label{rem:curvaturesign}
One may consider that $\nabla$ is \emph{positively curved} at $x$ if $S(\nabla,g)$ is positive at $x$ and \emph{negatively curved} if it is negative.
The sign of $S(\nabla,g)$ is the sign of $\psi$ and is thus locally constant (recall that the potential $\psi$ cannot vanish on the domain of definition of $\nabla$).
\end{rem}

The connection preserves any two-dimensional affine plane $P$ in an affine chart (see Corollary~\ref{cor:plane}), and if the restriction of $g$ to $\overrightarrow P$ is nondegenerate, then the restriction of $\nabla$ to $P$ satisfies $S(\nabla^P,g) = 2/\psi$, so has the same sign.
In the case of a general pseudo-Riemannian metric, in dimension $d=2$, the classical scalar curvature is $2$ times the Gaussian curvature. 
From this, it makes sense for $d=2$ to consider $1/\psi$ as the Gaussian curvature of $\nabla$ relative to $g$ and denote it $G(\nabla,g)$.

\section{Canonical models for pseudo-Euclidean manifolds}\label{sec:LocalModels}

In this section we work in local charts of pseudo-Euclidean manifolds. Therefore, we will consider a pseudo-Euclidean affine space $E$ of dimension $d$ with associated quadratic form $q$ of signature $(p,d-p)$ and with a hybrid connection $\nabla$ defined on a connected open subset of $E$. We only consider the case $d\geq 2$.

By Theorem~\ref{thm:MAIN}, the hybrid potential of $\nabla$ must satisfy $H(\psi)=q$, so it is of the form $q+A$ where $A$ is affine. Since $q$ is non-degenerate, an appropriate translation of $E$ transforms it into the form $q+\lambda$ for some $\lambda\in\R$.

\subsection{Canonical models of hybrid connections}\label{sub:canonical}

\begin{figure}
\begin{tikzpicture}
\node at (0,0) {\includegraphics[scale=0.4]{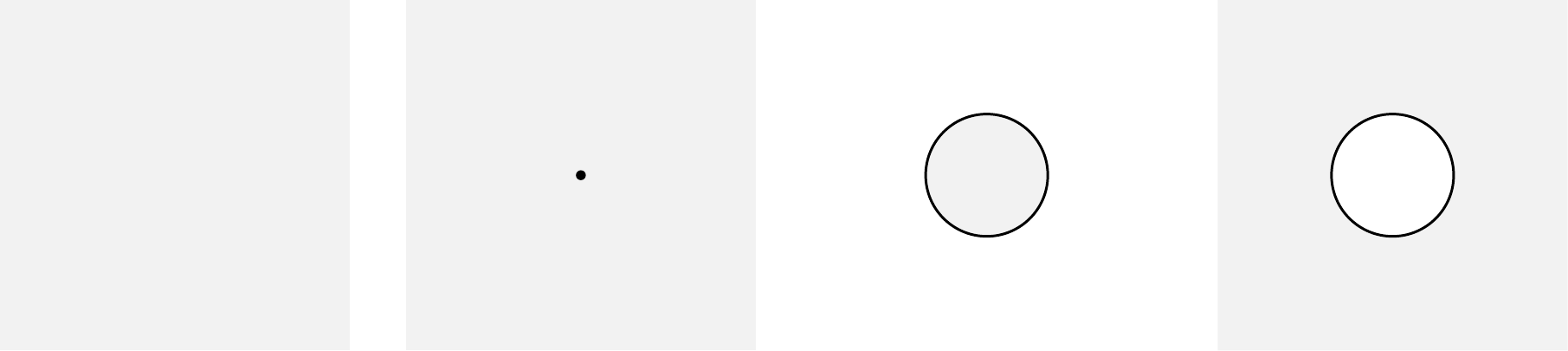}};
\node at (-4.7,-1.8) {$\mathcal{S}_{1}^{+}(2,0)$};
\node at (-1.55,-1.8) {$\mathcal{S}_{0}^{+}(2,0)$};
\node at (1.55,-1.8) {$\mathcal{S}_{-1}^{-}(2,0)$};
\node at (4.7,-1.8) {$\mathcal{S}_{-1}^{+}(2,0)$};
\end{tikzpicture}
\caption{Schematic presentation of the four models of Section~\ref{sub:canonical} for signature $(2,0)$.
The gray set shows the intersection of a square centered on $0$ and of the domain of the chart.
The black set its boundary in $\RR^2$.}
\label{fig:models:2:0}
\end{figure}
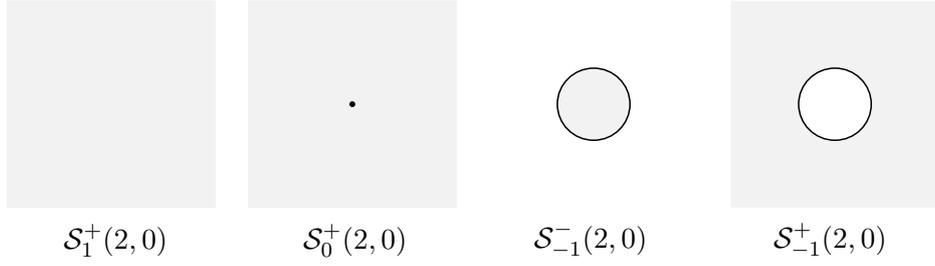

For each possible signature, we define here \textit{canonical models} for pseudo-Euclidean charts, corresponding to hybrid potentials $\psi(x)=q(x)+\lambda$ where $\lambda\in\R$.
\par
The connection is only defined on the set of points $x\in E^{p,d-p}$ where $\psi(x)\neq 0$. For a given potential $\psi$, we may have more than one connected component of the domain where $\psi\neq 0$.

\begin{defn}\label{defn:U+U-}
We assume $d\geq 1$. Let $E=E^{p,d-p}$.
We denote
\[U^{+} = \big\{x\in E\,;\, q(x)+\lambda>0\big\},\qquad
U^{-} = \big\{x\in E\,;\, q(x)+\lambda<0\big\}
.\]    
\end{defn}

We recall without proof the following basic facts in the next lemmas, that the reader can skip on first reading. Concerning non-emptyness:

\begin{samepage}
\begin{lem}\label{lem:bpe}
$U^{+}$ and $U^{-}$ satisfy the following statements:
\begin{itemize}
\item $U^{+}$ is empty if and only if $p=d$ and $\lambda\leq 0$.
\item $U^{-}$ is empty if and only if $p=0$ and $\lambda\geq 0$.
\item $U^+$ contains $0$ if and only if $\lambda>0$.
\item $U^-$ contains $0$ if and only if $\lambda<0$.
\item If $U^\pm$ contains $0$ then it deformation retracts onto $0$. Otherwise it deformation retracts to a subset homeomorphic to the unit sphere $S^{n-1}$ of the Euclidean space of dimension $n=p$ for $U^+$ or $n=p-d$ for $U^-$.
\end{itemize}
\end{lem}
\end{samepage}

Concerning connectedness:

\begin{lem}\label{lem:BPE1}
The following facts hold:
\begin{itemize}
\item $U^{+}$ is disconnected if and only if $p=1$ and $\lambda\leq 0$. In this case it consists of two components $U^{+}_L$ and $U^{+}_R$ defined by $x_1>0$ and $x_1<0$.
\item $U^{-}$ is disconnected if and only if $d-p=1$ and $\lambda\geq 0$. In this case it consists of two components $U^{-}_L$ and $U^{-}_R$ defined by $x_d>0$ and $x_d<0$.
\item Both $U^+$ and $U^-$ are disconnected if and only if the signature is $(1,1)$ and $\lambda=0$. 
\item When connected, $U^+$ is simply connected unless $p=2$ and $\lambda\leq 0$.
When disconnected, its two components are simply connected.
\item When connected, $U^-$ is simply connected unless $d-p=2$ and $\lambda\geq 0$.
When disconnected, its two components are simply connected.
\end{itemize}
\end{lem}

Concerning boundedness:
\begin{lem}\label{lem:bdd}
With the convention that the empty set is bounded:
\begin{itemize}
\item $U^{+}$ is bounded if and only if $p=0$.
\item $U^{-}$ is bounded if and only if $p=d$.
\item $U^+ = E$ if and only if $p=d$ and $\lambda>0$.
\item $U^- = E$ if and only if $p=0$ and $\lambda<0$.
\end{itemize}
\end{lem}

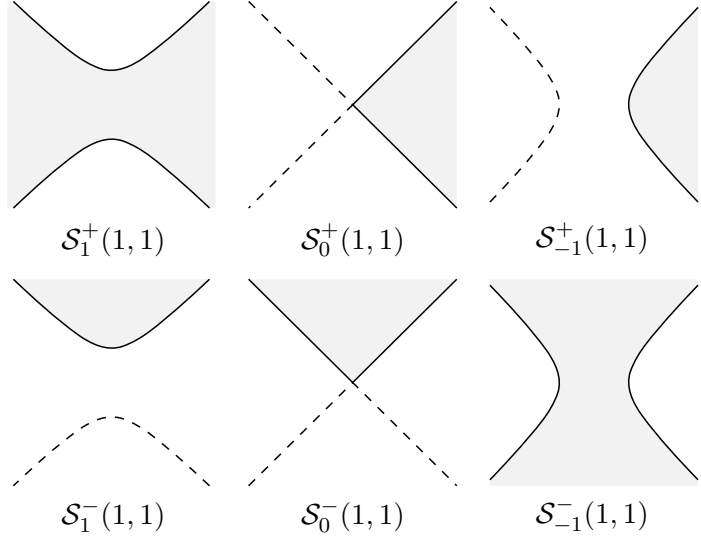
\begin{figure}[htbp]
\begin{tikzpicture}
\node at (0,-1.8) {\includegraphics[scale=0.4]{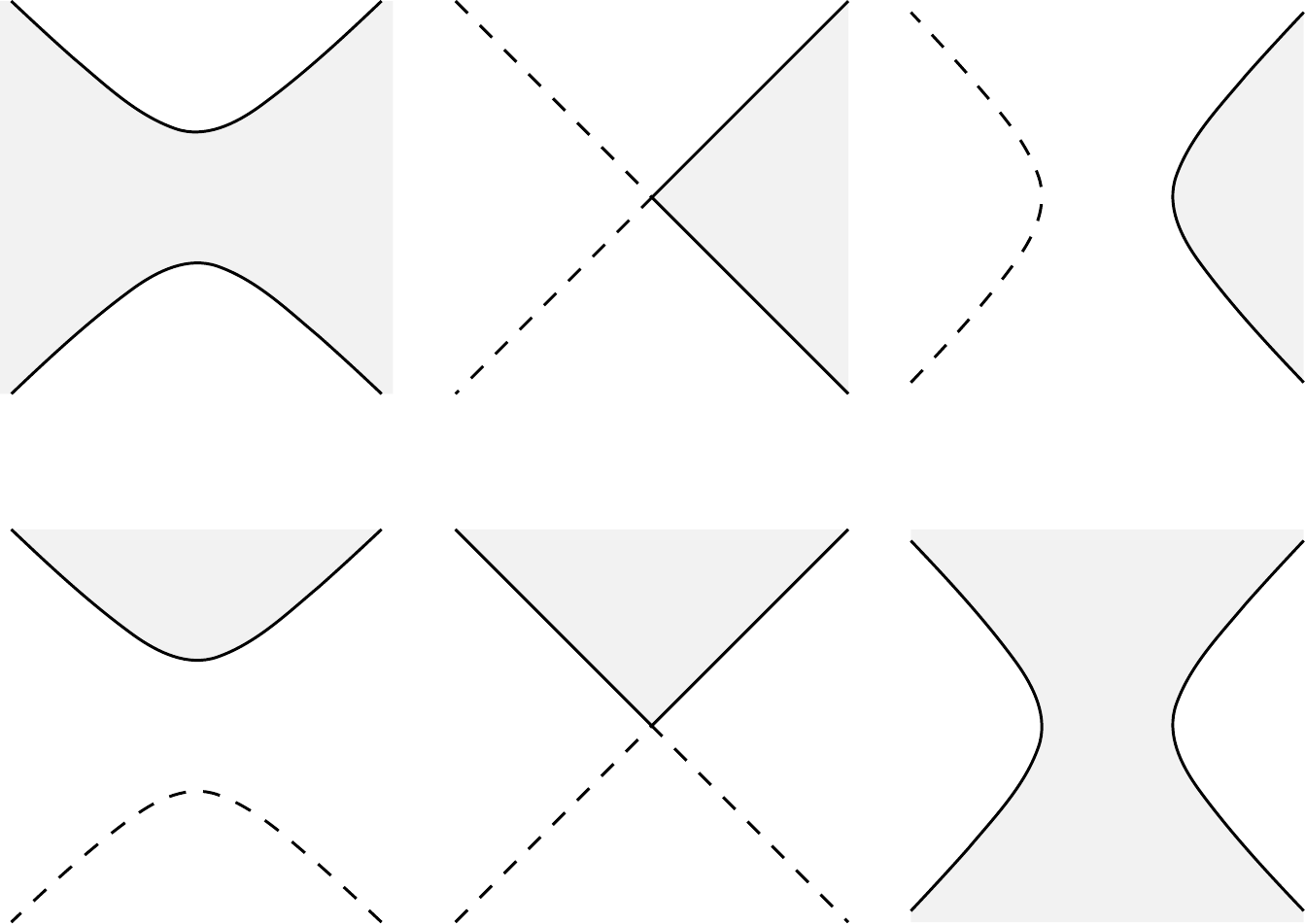}};
\node at (-3.15,-1.75) {$\mathcal{S}_{1}^{+}(1,1)$};
\node at (0,-1.75) {$\mathcal{S}_{0}^{+}(1,1)$};
\node at (3.15,-1.75) {$\mathcal{S}_{-1}^{+}(1,1)$};
\node at (-3.15,-5.4) {$\mathcal{S}_{1}^{-}(1,1)$};
\node at (0,-5.4) {$\mathcal{S}_{0}^{-}(1,1)$};
\node at (3.15,-5.4) {$\mathcal{S}_{-1}^{-}(1,1)$};
\end{tikzpicture}
\caption{Schematic presentation of the six models of Section~\ref{sub:canonical} for signature $(1,1)$ in $E^{1,1}$ (endowed with $q(x)={x_1}^2-{x_2}^2$).
The gray set shows the intersection of a square centered on $0$ and of the domain of the chart.
The black lines its boundary in $\RR^2$.
The dashed lines the remaining part of the locus of $\psi=0$.}
\label{fig:models:1:1}
\end{figure}


\begin{defn}[Canonical models]\label{def:models}
In the pseudo-Euclidean space $E^{p,d-p}$ endowed with the canonical form $q(x) = x_1^2+\cdots+x_p^2-x_{p+1}^2-\cdots-x_d^2$, let $\nabla$ be the hybrid connection of hybrid potential $\psi(x)=q(x)+\lambda$, $\lambda\in\R$.
We denote \[\mathcal{S}_{\lambda}^{+}(p,d-p)\] the restriction of $(E,\nabla)$ to the set of $x\in E$ such that $\psi(x)>0$ (i.e.\ the set denoted $U^{+}$ in the above lemmas), that we further restrict to $x_1>0$ if $p=1$ and $\lambda\leq 0$.
\par
Similarly, we denote
\[\mathcal{S}_{\lambda}^{-}(p,d-p)\]
the model restricted to $\psi<0$, and restricted furthermore to $x_d>0$ if $d-p=1$ and $\lambda\geq 0$.
\end{defn}

In particular, the canonical models are defined on connected components of the set $\big\{x\in E\,\vert \,\psi\neq 0\big\}$.

For signatures of the form $(d,0)$ and constants $\lambda \geq 0$, the set underlying $\mathcal{S}_{\lambda}^{-}(d,0)$ is empty, and we exclude this situation from the collection of models.
Similarly, for signatures of the form $(0,d)$ and constants $\lambda \leq 0$, the set underlying $\mathcal{S}_{\lambda}^{+}(0,d)$ is empty, and we also exclude this from the collection of models.
In all the other cases the model is non-empty and connected.
The models of signature $(2,0)$ are illustrated in Figure~\ref{fig:models:2:0} and for signature $(1,1)$ in Figure~\ref{fig:models:1:1}.

The function $\psi$ is invariant by the group $\on{O}(E)$ of the isometries of $E$ fixing $0$ (a.k.a. generalized orthogonal group).
Note that the models $\mathcal{S}_{0}^{\pm}(p,d-p)$ are invariant by the rescalings $x \mapsto rx$ with $r>0$.

\begin{rem}
There are correspondences between some models. For instance, there is a linear bijection between $\mathcal{S}_{\lambda}^{+}(p,d-p)$ and $\mathcal{S}_{-\lambda}^{-}(d-p,p)$ that multiplies $q$ by $-1$.
It is not a similarity.
\end{rem}

\begin{rem}\label{rem:ssng}
We introduced the scalar curvature of $\nabla$ relative to $g$ in Section~\ref{sub:ScalarCurvature} and denoted it $S(\nabla,g)$. We determined in Remark~\ref{rem:curvaturesign} that its sign is the sign of $\psi$. We immediately get that, on $\cal S^\nu_\lambda(p,d-p)$, it is $\nu$.
\end{rem}

\subsection{Slicing canonical models and the sign of curvature}\label{sub:chop}

If we consider a canonical model and cut it by an affine subspace $F$ of dimension at least two, intersecting the domain $U$ of the model, and such that the restriction of $q$ to $\vec F$ is nondegenerate, then one obtains a model isomorphic to one of the canonical models.

Indeed, by Corollary~\ref{cor:aff:restr} the connection has a restriction $\nabla^F$ to $U\cap F$ and the restriction of $\psi = q+\lambda$ to $F$ is a logarithmic potential for $\nabla^F$.
By nondegeneracy, $U\cap F$ and $\psi|_F$ can be mapped to a canonical model by some similarity/isometry.

Interestingly, all canonical models of dimension $d$ and arbitrary signature can be obtained as slices of a single pseudo-Euclidean space $E^{d+1,d+1}$ with the connection $\nabla$ associated with $\psi = q + 0$ and $U$ chosen as one of the two models it contains: $U=\cal S^+_0(d+1,d+1)$ and $U=\cal S^-_0(d+1,d+1)$.

We saw in Section~\ref{sub:ScalarCurvature} that the sign of the scalar curvature of $\nabla$ relative to $g$, $S(\nabla,g)$, is the same on the model restricted to an affine subspace $F$.
This sign is given by $\nu$ since $S(\nabla,g) = \frac{d(d-1)}{\psi}$.
On any two-dimensional affine plane $P$ through a point $x$ and such that the restriction of $q$ to $\ora P$ is nondegenerate, the interpretation of the sign of scalar curvature given in Section~\ref{sub:ScalarCurvature} can be rephrased as follows: an infinitesimal loop based on $x$ is an infinitesimal isometry pushing $v$ in the same sense as the loop if $\nu q(v) >0$ and in the opposite sense if $\nu q(v) <0$.

\subsection{Parallel transport in canonical models}\label{sec:ParallelTransport}

In any canonical model, we have $d\psi(x) = dq(x) = 2(x\cdot dx)$ where $x\cdot y$ denotes the associated scalar product of two vectors $x$ and $y$ with respect to the quadratic form $q$.
\par
Following Corollary~\ref{cor:expression} and Proposition~\ref{prop:Ricci}, at each point $x$ of a canonical model, the connection $\nabla$ is given by the following general formula:
\begin{equation}\label{eq:GF}
\Gamma(u,v)=-2\frac{(x \cdot u)v+(x \cdot v)u}{q(x)+\lambda}
\end{equation}

The parallel transport of a vector $v$ along a path $\gamma(t)$ is then characterized by the following differential equation:
\begin{equation}\label{eq:PT}
\dot{v}=-\Gamma(\dot{\gamma},v)=
\frac{2}{q(\gamma)+\lambda}\big((\gamma \cdot \dot{\gamma})v+(\gamma \cdot v)\dot{\gamma}\big)
\end{equation}

In this article, we call a vector $v$ a \emph{unit vector} if $q(v)\in\{-1,1\}$, and we call a family of unit vectors \emph{orthonormal} if they are pairwise orthogonal.

\subsection{The equipotential and radial foliations}\label{sub:foliation}

The \emph{light cone} is the set of $x\in E$ such that $q(x)=0$.
The \textit{radial foliation} is the foliation by straight half lines from $0$ (excluded).
In any canonical model $\mathcal{S}_{\lambda}^{\pm}(p,d-p)$, the level sets of the potential $\psi$ form the \textit{equipotential foliation}, whose leaves are quadric hypersurfaces, which are smooth, except possibly at $x=0$ if this point belongs to the model.
The equipotential foliation is preserved by the action of the generalized orthogonal group $\on{O}(p,d-p)$.
\par
An arc not contained in the light cone is a leaf of the radial foliation if and only if at each point it is orthogonal to the leaves of the equipotential foliation. Outside the light cone, the radial foliation and the equipotential foliation intersect transversely. The light cone is an equipotential leaf and is the union of $\{0\}$ and a collection of rays.

\subsection{Parallel transport along rays}\label{sub:PTrays}

In this Section, we consider the particular case of the parallel transport along a ray that does not belong to the light cone, and we parameterize it as $\gamma(t)=t e_{r}$ where $t$ ranges over some interval $[t_0,t_1] \subset \mathbb{R}^*_{+}$ and $e_{r}$ is a vector with $q(e_r)=\pm 1$. 

Any vector $v$ can be written in the form $v_{r}e_{r}+v_{\theta}e_{\theta}$ where $e_{\theta}$ is an orthoradial vector (tangent to the equipotential foliation at $\gamma(t_0)$).
The parallel transport of $v$ along a segment of a ray remains in the plane generated by $e_{r}$ and $e_{\theta}$ (see Corollary~\ref{cor:plane}).
We denote $v=v(t) = v_r e_r + v_\theta e_\theta$ the transported vector with $v_r=v_r(t)$ and $v_\theta=v_\theta(t)$.

\medskip

If $v$ is a radial vector, then $v_\theta=0$ and Equation~\eqref{eq:PT} simplifies to:
$$
\dot{v_{r}}=
\frac{4t q(e_{r})}{\lambda+t^{2}q(e_{r})}v_{r}
\text{~and~}
\dot{v_{\theta}}=0
$$
The vector remains radial all along the path and the solutions are then $v_\theta(t)=0$ and $v_{r}(t)=k(\lambda+t^2q(e_r))^2 = k(\lambda+q(x))^{2}$ with $x=\gamma(t) = t e_r$, where $k=1/(\lambda+q(\gamma(t_0)))^2\in\RR^*$.
In particular, we have the preserved quantity
\[\frac{1}{\psi(\gamma(t_1))^2}v(t_1) = \frac{1}{\psi(\gamma(t_0))^2} v(t_0) 
.\]

\medskip

If $v$ is an orthoradial vector then $v_r=0$ and we have $x \cdot v =0$ and the equation simplifies to:
$$
\dot{v_\theta}=\frac{2tq(e_{r})}{\lambda+t^{2}q(e_{r})} v_\theta,\quad \dot v_r = 0
.
$$
The transported vector remains orthoradial all along the path and the solutions are $v_r(t)=0$ and $v_{\theta}(t) = k(\lambda+t^2 q(e_\theta)) =k(\lambda+q(x))$ where $k=1/\lambda+q(\gamma(t_0))$.
In particular, we have the preserved quantity
\[\frac{1}{\psi(\gamma(t_1))} v(t_1) = \frac{1}{\psi(\gamma(t_0))} v(t_0) 
.\]

\medskip

Note in particular, that the exponent of the $\psi(\gamma(t))$ factors is different for radial and orthoradial vectors.

\begin{prop}\label{prop:PTrays}
Consider a ray $\cal R$ not contained in the light cone and let $e_r$ be a unit vector in $\cal R$.
In an orthonormal basis where $e_r$ is the first vector, the parallel transport along a path $\gamma$ contained in $\cal R$ corresponds to a diagonal matrix:
\[
  \on{Diag}(\Lambda^2,\Lambda,\cdots,\Lambda)
\]
where $\Lambda = \frac{\psi(\gamma(t_1))}{\psi(\gamma(t_0))}$. The radial and every orthoradial direction are preserved. In particular, the tangent hyperplane to the equipotential foliation is preserved.
\end{prop}

\medskip
\noindent\textsf{Parallel transport along rays in the light cone.} The domain of a canonical model is the locus where $\lambda+q(x) \neq 0$.
In particular, if $\lambda \neq 0$ and the sign of $\lambda$ is $\nu$, a point $x$ in the model may belong to the light cone based on $0$, i.e.\ the locus where $\psi=\lambda$.
The computations performed above are still valid, and since $\psi$ is constant (and equal to $\lambda$) on the ray, they yield that $v(t)$ remains constant provided $v(0)$ is radial or orthoradial.

The only difference is that the radial and orthoradial subspaces no longer generate the whole vector space: actually the radial direction is now contained in the orthoradial hyperplane.

To handle the general case, let us write $\gamma(t) = t e_r$ with $e_r\cdot e_r =q(e_r) =0$, so $\dot \gamma = e_r$ and $\gamma \cdot \dot \gamma = 0$.
From Equation~\eqref{eq:PT}, for any initial vector $v(t_0)$, the transport equation reads
\[\dot{v} = 2t\frac{e_r \cdot v}{\lambda}e_r.\]
It follows that $v(t) = v(t_0) + f(t) e_r$.
Knowing this, the transport equation is equivalent to
\[\dot{f} = 2t\frac{e_r \cdot v(t_0)}{\lambda}.\]
Hence
$$
v(t) = v(t_0) + \frac{e_r \cdot v(t_0)}{\lambda}(t^{2}-t_0^{2})e_{r}
.$$

By properties of the parallel transport of connections, using a non-uniformly parameterized ray $\gamma(t) = f(t) e_r$ does not basically change the result. Here is a summary:

\begin{cor}\label{cor:PTlight}
For any path $\gamma$ contained in a radial leaf contained in the light cone:
\begin{itemize}
\item The parallel transport of a vector tangent to the light cone along $\gamma$ is trivial (i.e. constant).
\item Writing $\gamma(t) = f(t) e_r$ with $q(e_r)=0$, the parallel transport of a general vector $v$ along $\gamma$ satisfies
\[ v(t) = v(t_0) + \frac{e_r \cdot v(t_0)}{\lambda}(f(t)^{2}-f(t_0)^{2})e_{r} 
.\]
\end{itemize}
\end{cor}

\subsection{Parallel transport along equipotential arcs}\label{sub:PTcirc}

Outside the light cone, equipotential leaves are quadric hypersurfaces of $E$.
We will not try to compute the parallel transport along any path contained in such leaves, but rather focus on planar sections of these leaves. By \emph{plane} we mean $2$-dimensional vector (or affine) subspaces of $E$.

\begin{defn}
In any canonical model, we call an \textit{equipotential arc} a non-empty portion of the intersection of an equipotential leaf different from the light cone with a plane $\mathcal{P}$ containing $0$ and whose induced pseudo-Euclidean structure is nondegenerate.
\end{defn}

Every point not in the light cone belongs to at least one equipotential arc, as follows from the following (slightly more general) lemma:

\begin{lem}
If $d\geq 2$, then for any $x\neq 0$  there exists a plane through $x$ and $0$ on which $q$ is nondegenerate.
\end{lem}
\begin{proof}
If $q(x)\neq 0$ then $\R x$ is transverse to its orthogonal complement $F$, on which $q$ cannot be zero, for otherwise the scalar product would be zero too on $F$, and hence any $y\in (\R x)^\perp$ would be in the kernel of $q$, contradicting nondegeneracy of $q$ on $E$; taking any $y\in F$ such that $q(y)\neq 0$, we get that the matrix of $q$ in the basis $(x,y)$ of $\R x+\R y$ is diagonal with non-vanishing coefficients.
If $q(x)=0$ then since $x\notin \{0\}=\on{Ker} q$, there exists $y\in E$ such that $x\cdot y\neq 0$. Then the determinant of the matrix of $q$ in the basis $(x,y)$ of $\R x+\R y$ is equal to $-(x\cdot y)^2 \neq 0$.
\end{proof}

On an equipotential arc $\gamma$ of the model $\cal S^\nu_\lambda(p,d-p)$ with $\nu\in\{-,+\}$ we have $q(\gamma)=\on{cst}\neq 0$ and we let $r>0$ be such that\footnote{Where $\nu=-$ is replaced by $-1$ and $\nu=+$ by $1$.}
\[q(\gamma) = \nu r^2.\]
Depending on whether the signature of the pseudo-Euclidean structure of the plane is $(2,0)$, $(0,2)$ or $(1,1)$, an equipotential arc $\gamma$ is parameterized by a usual angle or a hyperbolic angle $\theta$, characterized as follows: $|q(\dot{\gamma}(\theta))|=|q(\gamma(\theta))| = r^2$.
Since $\gamma$ is an equipotential arc, the sign of $q(\dot\gamma)$ is determined in advance: $q(\dot{\gamma}(\theta))=\eps\,q(\gamma(\theta))$ with $\eps = -1$ in signature $(1,1)$ and $\eps = 1$ otherwise.
Yet, this still allows for two opposite values of $\dot\gamma(\theta)$.
To fix one, we choose one of the two possible orientations on the curve $\cal L$, where $\cal L$ is the component containing $\gamma$ of the equipotential leaf in plane $\cal P$. 
Then we impose $\dot\gamma(\theta)$ to be in the positive sense of this orientation.

\medskip

By the nondegeneracy assumption and because we are not on the light cone, there is a moving orthonormal frame $(e_{r},e_{\theta})$ of $\mathcal{P}$ so that
\begin{itemize}
\item $\gamma(\theta) = r e_r$;
\item $\dot\gamma(\theta) = r e_\theta$.
\end{itemize}
As is customary for instance in physics, the symbol $r$ plays two different roles, as does $\theta$. As indices of $e$ here or $v$ below, they are purely symbolic and allow us to distinguish one object from another.
We have
\begin{itemize}
\item $\dot{e_r} = e_\theta$;
\item $\dot{e_\theta} = -\epsilon e_r$;
\end{itemize}
where\footnote{We could have taken the opposite convention for $\eps$ so as to have the simpler expression $\dot{e_\theta} = \eps e_r$ but our choice makes most other expressions simpler.}
\begin{itemize}
\item $\epsilon = 1$ in signatures $(2,0)$ and $(0,2)$;
\item $\epsilon = -1$ in signature $(1,1)$.
\end{itemize}

\medskip

We first observe the following fact. 
By the nondegeneracy assumption, $\cal P$ is transverse to its orthogonal space in $E^{p,d-p}$.
For any vector $v$ that is orthogonal to $\mathcal{P}$, the parallel transport along $\gamma$ of $v$ is trivial.
Indeed, the vector $v$, the position vector $\gamma$ and the speed vector $\dot{\gamma}$ are orthogonal to each other. Equation~\ref{eq:PT} thus reduces to $\dot{v}=0$. In other words, the parallel transport along an equipotential arc contained in a plane $\mathcal{P}$ preserves any vector orthogonal to $\mathcal{P}$.

\medskip

Now we consider a vector $v$ of the form $v_{r}e_{r}+v_{\theta}e_{\theta}$. The parallel transport of $v$ along $\gamma$ always remains in $\mathcal{P}$ (by Corollary~\ref{cor:plane}).
Equation~\eqref{eq:PT} for the parallel transport of $v$ along  $\gamma$ reduces to:
$$\dot{v}=\frac{2 (v \cdot \gamma)\dot{\gamma}}{\lambda+q(\gamma)}=\frac{2q(\gamma)}{\lambda+q(\gamma)}v_{r}e_{\theta}
.$$
Recall that $q(\gamma)$ is constant.
Since $\dot{v}=(\dot{v_{r}}-\epsilon v_{\theta})e_{r}+(v_{r}+\dot{v_{\theta}})e_{\theta}$, we obtain a system of two differential equations:
\begin{itemize}
    \item $\dot{v_{r}}-\epsilon v_{\theta}=0$;
    \item $v_{r}+\dot{v_{\theta}}
=\frac{2q(\gamma)}{\lambda+q(\gamma)}v_{r}
    $.
\end{itemize}
The system is then equivalent to:
\begin{itemize}
    \item $\epsilon v_{\theta}= \dot{v_{r}}$;
    \item $\Ddot{v_{r}}=\epsilon\frac{q(\gamma)-\lambda}{q(\gamma)+\lambda}v_{r}$.
\end{itemize}
This is the \emph{harmonic oscillator} equation with position variable $v_r$ and speed variable $\epsilon v_\theta$.

For any equipotential leaf (where $q$ is constant), we introduce the \textit{characteristic frequency}
\[\omega = \sqrt{\epsilon \frac{q(\gamma)-\lambda}{q(\gamma)+\lambda}} = \sqrt{\epsilon \frac{1-\lambda'}{1+\lambda'}}\in\CC\]
with $\lambda' := \lambda / q(\gamma) = \nu\lambda/r^2$.
Recall that the connection is not defined on the locus where $q+\lambda$ vanishes, so $q(\gamma)+\lambda\neq 0$ in our case.

The quantity $\omega$, defined up to a sign, is either real or purely imaginary (we can fix $\omega$ by imposing that $\Re \omega$ and $\Im \omega$ be non-negative):
\begin{itemize}
    \item if $\lambda =0$, then $\omega = 1$ in signatures $(2,0)$ or $(0,2)$ and $\omega = i$ in signature $(1,1)$;
    \item if $q(\gamma)=\lambda$ then $\omega=0$; 
    \item if $|q(\gamma)| > |\lambda|$, then $\omega$ is a real number in signatures $(2,0)$ or $(0,2)$, purely imaginary in signature $(1,1)$;
    \item if $|q(\gamma)| < |\lambda|$, then $\omega$ is a purely imaginary number in signatures $(2,0)$ or $(0,2)$, real in signature $(1,1)$.
\end{itemize}

In the case $\omega \in \mathbb{R}$, solving the harmonic oscillator equation provides the following solution:
\[\begin{pmatrix}v_{r}(t) \\ v_{\theta}(t)\end{pmatrix}
=
\begin{pmatrix}
\cosh \omega t & \frac{\epsilon}{\omega}\sinh \omega t
\\
\epsilon\omega \sinh \omega t & \cosh \omega t
\end{pmatrix}
\cdot
\begin{pmatrix}
v_{r}(0)
\\
v_{\theta}(0)
\end{pmatrix}
\]

If $\omega$ is purely imaginary, we consider $\alpha=-i \omega$. Then $\alpha\in\RR$ and we obtain the following solution:
\[\begin{pmatrix}v_{r}(t) \\ v_{\theta}(t)\end{pmatrix}
=
\begin{pmatrix}
\cos \alpha t & \frac{\epsilon}{\alpha}\sin \alpha t
\\
-\epsilon \alpha \sin \alpha t & \cos \alpha t
\end{pmatrix}
\cdot
\begin{pmatrix}
v_{r}(0)
\\
v_{\theta}(0)
\end{pmatrix}
\]

Actually the two matrices coincide in all cases using complex versions of the trigonometric and hyperbolic functions.
Note also that choosing the opposite square root for $\omega$ does not change the matrices given by the formulae above, as expected since the parallel transport does not depend on the choice of $\omega$.

Do not forget that the matrices above are expressed in a polar coordinate system.
Unless we are in signature $(2,0)$ or $(0,2)$ and perform a full turn around the (circular) leaf, for the vector space structure of $\cal P$, this matrix is expressed in bases in the domain and range that are \emph{different bases} of the plane $\cal P$.
To recover its expression in a single basis of $\cal P$, one has to multiply it on each side with appropriate transition matrices.

Remarkably, if $\omega=0$ then the coordinates $v_r$ and $v_\theta$ are constant, i.e.\ $v$ is constant in the moving frame, but as we just saw, $v$ is not constant in a basis of $E$.
An equipotential leaf with $\omega=0$ exists on the models $\cal S^{+}_{1}(2,0)$, $\cal S^{-}_{-1}(0,2)$, $\cal S^{+}_{1}(1,1)$ and $\cal S^{-}_{-1}(1,1)$ but not on the other 10 canonical models in dimension 2.
In signatures $(2,0)$ and $(0,2)$ these are the models whose domain is the whole plane.
In signature $(1,1)$ these are the models whose domain contains $0$.

\medskip

\noindent\textsf{Performing a full turn along a circular leaf.} In the case of signature $(2,0)$, by taking $t=2\pi$, we obtain the holonomy matrix corresponding to a loop, more precisely a full turn along the circle of Euclidean radius $r$.

Remarkably, when $\lambda=0$, this holonomy matrix is independent of the leaf and is
\[
\begin{pmatrix}
\cosh 2\pi & \sinh 2\pi
\\
\sinh 2\pi & \cosh 2\pi
\end{pmatrix}
.\]
The numerical values are $\cosh 2\pi \approx 267.7467615$ and $\sinh 2\pi \approx 267.7448940$.

Still in the case $\lambda=0$ in signature $(2,0)$, note that the parallel transports along closed or non-closed arcs of circles are isometries for the pseudo-Riemannian structure $v_r^2-v_\theta^2 = (dr)^2-r^2(d\theta)^2$, whereas the Riemannian structure induced by the Euclidean norm is given by $v_r^2+v_\theta^2 = (dr)^2+r^2(d\theta)^2$.

In the case $\lambda=-1$ and $\nu=-1$, for which the domain of $\cal S$ is the Euclidean unit disk, if we take the circle of Euclidean radius $r$, we get the holonomy matrix
\[
M = \begin{pmatrix}
\cos 2\pi \alpha & \frac{1}{\alpha}\sin 2\pi \alpha 
\\
- \alpha \sin 2\pi\alpha & \cos 2\pi\alpha 
\end{pmatrix}
\]
with 
\[\alpha = \sqrt{\frac{1+r^2}{1-r^2}}
.\]
When $r\longrightarrow 0$ we get the expansion
\[
M = \on{I} - 2\pi r^2\begin{pmatrix}
0 & -1
\\
1 & 0
\end{pmatrix} + O(r^3)
.
\]
In other words, a transported vector asymptotically turns by $-2\pi r^2$.
This is coherent with the value $-2$, at $x=0$, of the Gaussian curvature of $\nabla$ relative to the flat metric $g$, which we found is equal to $2/\psi$. Here $\psi(0)=-1$ and the Euclidean area of the disk of Euclidean radius $r$ is $\pi r^2$.

\medskip

\noindent\textsf{Model $\cal S^{\pm}_0(1,1)$.} For $\lambda=0$ in signature $(1,1)$ the holonomy matrices for paths along leaves (different from the light cone) are the matrices
\[
\begin{pmatrix}
\cos \theta & -\sin \theta
\\
\sin \theta & \phantom{-}\cos \theta
\end{pmatrix}
\]
where $\theta$ measures the difference of hyperbolic angles along the leaf. The sign convention is with respect to the initially chosen orientation on the (one-dimensional) leaf $\cal L$.
These are isometries for the Riemannian structure $v_r^2+v_\theta^2 = (dr)^2+r^2(d\theta)^2$, whereas the pseudo-Riemannian structure induced by the pseudo-Euclidean norm is given by $\nu\times(v_r^2-v_\theta^2) = \nu\times((dr)^2-r^2(d\theta)^2)$.

\subsection{Local holonomy}\label{sub:LocHol}

The connection of a canonical model is incompressible, so the holonomy of contractible loops is contained in the group $\on{SL}(d,\mathbb{R})$.
We prove in this section that the \textit{local holonomy groups} $\on{Hol}(x,V)$ (generated by the holonomy of loops contained in a neighborhood $V$ of $x$) always coincide with $\on{SL}(d,\mathbb{R})$.

\begin{prop}\label{prop:LocHol}
In any canonical model $\cal S=\mathcal{S}_{\lambda}^{\pm}(1,1)$ or $\cal S=\mathcal{S}_{\lambda}^{\pm}(2,0)$, for any open subset $V$ of $\cal S$ and any $x\in V$, the holonomy group $\on{Hol}(x,V)$ of the connection is equal to $\on{SL}(2,\mathbb{R})$.
\end{prop}

\begin{proof}
We already know that the holonomy locally preserves some volume form and is therefore contained in $\on{SL}(2,\mathbb{R})$.
The parallel transport locally conjugates the holonomy groups at different points of $V$.
It is hence enough to consider the case where $x$ is not on the light cone (whose interior is empty).
By Theorem~4.2 of \cite{KN}, $\on{Hol}(x,V)$ is a Lie subgroup of $\on{SL}(2,\mathbb{R})$.
Its Lie algebra $\cal L$ contains the infinitesimal holonomies $M(i,j)$ mentioned in Section~\ref{sec:AffCon}, whose matrices satisfy $M(i,j)_{kl} = -R_{ijk}^l$.
By Equation~\ref{eq:Rpsi}, $M(1,2)$ has matrix $\frac{1}{\psi} J_\sigma$ with $J_\sigma = \begin{pmatrix} 0 & -\epsilon_2  \\ \epsilon_1 & 0 \end{pmatrix}$.
In particular, $J_\sigma\in\cal L$.
Consider now the parallel transport along the ray through $x$ and its expression in an adapted basis $(e_r,e_\theta)$.
This new basis is also orthonormal and the change of basis turns $J_\sigma$ into $P J_\sigma P^{-1}$ which equals either $J_\sigma$ or $-J_\sigma$.
Let $\tau>1$ be close enough to $1$.
According to \ref{sub:PTrays}, the parallel transport from $\tau x$ to $x$ has, in the basis $(e_r,e_\theta)$, matrix $D=\on{Diag}(\Lambda^2,\Lambda)$ where $\Lambda = \frac{\psi(x)}{\psi(\tau x)}$. Note that $\Lambda$ is close to $1$, but different.
If we transport the infinitesimal holonomy element $J_\sigma$ from $\tau x$ to $x$ we obtain $D J_\sigma D^{-1} = \begin{pmatrix} 0 & -\Lambda\epsilon_2  \\ \Lambda^{-1}\epsilon_1 & 0 \end{pmatrix}$.
Differentiating as $\tau$ tends to $1$ we get that $K_\sigma\in P\cal L P^{-1}$ where $K_\sigma = \begin{pmatrix} 0 & \epsilon_2  \\ \epsilon_1 & 0 \end{pmatrix}$.
One computes
$[J_\sigma,K_\sigma] = -2\epsilon_1\epsilon_2 L_\sigma$ with $L_\sigma = \begin{pmatrix} 1 & 0 \\ 0 & -1 \end{pmatrix}$, hence $ L_\sigma \in P\cal L P^{-1}$.
This proves the proposition, since $(J_\sigma,K_\sigma,L_\sigma)$ is a basis of the Lie algebra of $\on{SL}(2,\mathbb{R})$.
\end{proof}

\begin{thm}\label{thm:LocHol}
In any canonical model $\mathcal{S}_{\lambda}^{\pm}(p,d-p)$, for any open subset $V$ and any $x\in V$, the  holonomy group $\on{Hol}(x,V)$ of the connection is $\on{SL}(d,\mathbb{R})$.
\end{thm}

\begin{proof}
We already know that the holonomy locally preserves some volume form hence is contained in $\on{SL}(d,\mathbb{R})$.
Again, it is enough to compute the holonomy group in  neighborhood of a point $x$ that does not belong to the light cone.
Then the line generated by the vector $x-0$ is disjoint from its orthogonal hypersurface $H$ (the tangent space at $x$ to the foliation is $x+H$).
The form $q$ is nondegenerate on $H$ (otherwise it would be degenerate on $\R^d$) so we can choose an orthonormal basis $e_2,\dots,e_d$ of $H$, and complete it into an orthonormal basis $e_{1},\dots,e_{d}$ of $\R^d$ where $e_{1}$ is a radial vector.
We denote by $\mathcal{L}$ the Lie algebra that is tangent to the holonomy group $\on{Hol}(x,V)$ at the identity map.
We consider the holonomy of loops contained in a plane $\cal P$ through $x$ and directed by the vector plane $P$ generated by $e_{1}$ and some other element $e_{i}$ of the basis. 
Following Sections~\ref{sub:PTrays} and~\ref{sub:PTcirc}, these holonomies are the identity on $P^\perp$.
\newline
On $\cal P$, the parallel transport is equivalent to a nondegenerate two-dimensional model.
From this and Proposition~\ref{prop:LocHol}, we get that for any $2 \leq i \leq d$, $\mathcal{L}$ contains a copy of $\mathfrak{sl}(2,\RR)$.
More precisely, the following basis elements of $\mathfrak{sl}(2,\RR)$: $\begin{pmatrix} 1 & 0 \\ 0 & -1 \end{pmatrix}$, $\begin{pmatrix} 0 & 1 \\ 0 & 0 \end{pmatrix}$, $\begin{pmatrix} 0 & 0 \\ 1 & 0 \end{pmatrix}$ correspond respectively, in $\mathfrak{sl}(d,\R)$, to $e_{1} \otimes e_{1}^{\ast}-e_{i} \otimes e_{i}^{\ast}$, $e_{1} \otimes e_{i}^{\ast}$ and $e_{i} \otimes e_{1}^{\ast}$.
\newline
For $i \neq j$, the Lie bracket of $e_{j} \otimes e_{1}^{\ast}$  and $e_{1} \otimes e_{i}^{\ast}$ equals $e_{i} \otimes e_{j}^{\ast}$.
It follows that the Lie algebra $\mathcal{L}$ contains the basis $e_{1} \otimes e_{1}^{\ast}-e_{i} \otimes e_{i}^{\ast}$ ($i\neq 1$), $e_{i} \otimes e_{j}^{\ast}$ ($i\neq j$), of $\mathfrak{sl}(d,\mathbb{R})$.
The local holonomy group of the connection thus coincides with $\on{SL}(d,\mathbb{R})$.
\end{proof}

\section{Properties of hybrid connections on pseudo-Euclidean manifolds}\label{sec:Tale}


Recall that a \emph{conformal class} is an equivalence class of pseudo-Riemannian metric $g$ up to multiplication by a positive function $f$, i.e.\ $g\sim f g$.

\medskip

Any canonical model $\cal S^{\nu}_{\lambda}(p,q)$, $\nu \in\{-,+\}$ with its connection $\nabla$ takes place within a pseudo-Euclidean space $E^{p,d-p}$, endowed with a quadratic form $q(x)=\sum \eps_i\, {x_i}^2$, the associated scalar product $(x,y)\in E^2\mapsto x\cdot y\in \R$ and the associated flat pseudo-Riemannian tensor
\[g(dx,dy) = dx \cdot dy.\]
However, $g$ is not preserved by the connection $\nabla$.

On the opposite, $\nabla$ has an unambiguous\footnote{In the sense that all authors agree on the sign in its definition.} Ricci tensor
\[
\Ric=\frac{d-1}{\psi}g
\]
(see Corollary~\ref{cor:HybridRic}).
If $\nu<0$ then $\psi<0$, so the signature of the Ricci tensor is $(d-p,p)$ (the opposite signature).
If $\nu>0$ then $\psi>0$, so the signature of the Ricci tensor is $(p,d-p)$ (the same as for $g$).

\begin{rem}
In \cite{Ko} have been given many examples of incompressible connections with positive definite or negative definite Ricci tensors.
Our canonical models provide another family of examples when $(p,d-p) = (d,0)$ or $(0,d)$.
\end{rem}

As we will see in Section~\ref{sub:isochrone}, we can define another pseudo-Riemannian metric $h$ with remarkable properties: it is in the conformal class of $g$, the geodesics of $\nabla$ have constant speed with respect to $h$, and its expression is
\[ h=\frac{g}{\psi^{4}}
\]
where $\psi$ is the potential.

In the end, we find in the same conformal class three metrics:
\begin{enumerate}
    \item the flat pseudo-Euclidean metric $g$;
    \item the Ricci tensor if $\psi>0$ (or its opposite if $\psi<0$) $\Ric=\frac{d-1}{\psi}g$;
    \item the metric $h=\frac{g}{\psi^{4}}$.
\end{enumerate}

Formally, Propositions~\ref{prop:Ricci} and~\ref{prop:metr} form a proof of Theorem~\ref{thm:ThreeMetrics}. The remainder of this section is devoted to making explicit the properties of the metric $h=\frac{g}{\psi^{4}}$ and their consequences for the parameterization of geodesics, in particular the special case in which the hybrid connection is geodesically complete.

\subsection{Isochrone metric}\label{sub:isochrone}

\begin{defn}\label{def:isochrone}
A metric and a connection are \emph{isochrone} to each other if the geodesics of the connection have constant speed with respect to the metric.
\end{defn}

Note that geodesics for the metric and for the connection do not have to coincide.

\begin{prop}\label{prop:metr}
In any pseudo-Euclidean manifold $(M,D,g)$, the hybrid connection $\nabla$ defined by a potential $\psi$ is isochrone to
$$ h=\frac{\alpha^2}{\psi^{4}}g $$
for any $\alpha>0$.
\par
Moreover, these are the only isochrone metrics for $\nabla$ in the same conformal class as $g$ (by this we mean metrics of the form $\rho(x)g$ where $\rho$ is a positive function).
\par
For this reason, and because the hybrid potential $\psi$ is uniquely defined, we refer to $\frac{g}{\psi^{4}}$ as \emph{the isochrone metric} of $\nabla$.
\end{prop}

\begin{proof}
Since any geodesic in an affine chart is contained in a straight line $L\subset \R^d$, the geodesic equation $\ddot \gamma = -\Gamma(\dot \gamma,\dot\gamma)$ is easy to solve.
\par
From $\Gamma(\dot \gamma,\dot\gamma) = 2a(\dot\gamma)\dot\gamma = -\frac{2}{\psi} d\psi(\dot\gamma)\dot\gamma$, the geodesic equation reads
\[ \ddot\gamma = 2 \frac{\dot{\wideparen{\psi(\gamma)}}}{\psi(\gamma)}\dot\gamma
,\]
which is equivalent to $\frac{d}{dt}(\dot\gamma/\psi(\gamma)^2)=0$, i.e.\ 
\[ \frac{1}{\psi(\gamma)^2}\dot \gamma = \on{cst}
.\]
If one chooses an affine parameterization of $L$, $\phi : s\in\R\mapsto x\in L$, and denote $\psi_s = \psi\circ\phi$ (the subscript $s$ is purely symbolic) then we get the equation
$\dot s /\psi_s^2(s) = \on{cst}$.
This is equivalent to
\[F(s(t)) = \alpha t+\beta\]
for some $\alpha\in\R^*$ and $\beta\in\R$, 
where $F$ is an antiderivative of $1/\psi_s^2$. 
This can be rephrased as
\[\int_{s(t_0)}^{s(t_1)}\frac{ds}{\psi_s(s)^2} = \alpha (t_1-t_0);\]
or as follows: the change of variables $F\circ\phi^{-1}$ with $F' = 1/\psi_s^2$ trivializes the restriction of $\nabla$ to $L$.
\par
Then, for any pseudo-Euclidean metric $g$ on $M$ (i.e.\ with $g$ constant in affine charts), every geodesic of $\nabla$ has constant speed with respect to the pseudo-Riemannian metric $h$ such that $h(dx,dx) = \frac{g(dx,dx)}{\psi(x)^{4}}$ so $\nabla$ and $h$ are isochrone to each other. 
\par
Conversely any two isochrone and $g$-conformal metrics of the form $\rho g$ and $\tilde \rho g$ must have a locally constant ratio $\alpha=\tilde \rho/\rho$ on lines whose direction is not in the light cone.
Locally any point close enough to $x$ can be reached by following one or two straight segments whose directions are not in the light cone.
It follows that $\tilde \rho/\rho$ is locally constant.
\end{proof}

\begin{rem}\label{rem:lcgeo}
Recall that the geodesics of $\nabla$ are supported by straight lines in $E$.
If the geodesic direction $\dot \gamma$ is in the light cone of $q$, i.e.\ if $q(\dot \gamma)=0$, then the speed with respect to $g$, hence with respect to $h$, is automatically $0$ along $\gamma$, hence the isochrone property does not bring any information for these specific geodesics.
However, for any straight line $L$ of $E$ whose direction $\vec L$ is not contained in the light cone of $q$, and with non-empty intersection with the domain $U$ of the canonical model, and for any choice of Euclidean metric\footnote{We call a metric on $L$ Euclidean if it has constant coefficients in the coordinates $(x_{1},\dots,x_{d})$ on $E$.} $\tilde g$ on $L$, the geodesics supported by $L$ are isochrone with respect to the metric $\tilde g/\psi^4$.
By passing to the limit on the direction of $\vec L$, this also applies to the case when $\vec L$ is contained in the light cone of $q$.
\end{rem}

\begin{rem}
In general, the parallel transport under $\nabla$ does not consist in isometries of the isochrone metric $h = \frac{g}{\psi^{4}}$ of Proposition~\ref{prop:metr}.
In other words, the symmetric connection $\nabla$ does not coincide with the Levi-Civita connection of $h$.
\end{rem}

\begin{rem}
Consider a canonical model and the isochrone metric $h=\frac{g}{\psi^4}$ defined with respect to $\nabla$.
It would be tempting to define a distance adapted to $\nabla$ as follows:
given two points $a$ and $b$ let $d(a,b)$ be the distance between $a$ and $b$ with respect to the metric $h$ restricted to $[a,b]$.
However there are three problems:
\begin{itemize}
\item one has to assume that the straight segment $[a,b]$ is contained in the domain of the model, so not every pair of points can necessarily be joined, depending on the model;
\item the triangle inequality is not necessarily satisfied;\footnote{See Section~\ref{sub:ex:isoc} for counterexamples.}
\item in mixed signatures, the distance would be positive, zero or imaginary, depending on the sign of $q(b-a)$.
\end{itemize}    
\end{rem}

\subsection{More about the isochrone metric}\label{sub:more:isoc}

In this section we investigate properties of the metric $h=g/\psi^4$.
First, we would like to compute the scalar curvature $S(h)$ of $h$.
We recall that this is a $(0,0)$-tensor (that is, a function), i.e.\ it is independent of the chart.
If one chooses another isochrone metric by multiplying $h$ by $\alpha^2$, this multiplies $S$ by $\alpha^{-2}$.

To compute the scalar curvature we use the formula (see \cite{Be}, Chapter~1, Section~J, Theorem~1.159)
\[ S(e^{2f} g) = e^{-2f}\big(S(g)-2(d-1)\Delta_g f-(d-2)(d-1)g(df,df)\big)
\]
where $g$ is a general pseudo-Riemannian metric, $\Delta_g$ is the pseudo-Laplacian associated with $g$ and the musical isomorphism is used to define $g(df,df)$.
In the case of our flat $g=\sum_i \eps_i\, dx_i\otimes dx_i$, we have $S(g)=0$ and $\Delta_g = \sum_i \eps_i \frac{\partial^2}{(\partial x_i)^2}$. Then, letting $f=-2\log|\psi|$ and remembering that $\psi = q +\lambda$ we get:
\[ S(h) = 8(d-1)(d\lambda-(d-2)q)\psi^2
.\]
Remarkably, $S(h)$ vanishes on the boundary of the canonical model.

In dimension $d\geq 3$, the sign of $S(h)$ will change on the set defined by the equation $q(x)=\frac{d}{d-2}\lambda$, which is non-empty in the domain $U$ of $\cal S^\nu_\lambda(p,d-p)$ if and only if $\lambda\neq 0$, $\lambda$ and $\nu$ have the same sign and there is at least one $\eps_i$ with this sign.
In this case it is a smooth equipotential.
Moreover, $U$ contains $0$ and $S(h)$ has the sign of $\lambda$ at $x=0$.
If one replaces $\nu$ by $-\nu$, assuming that the domain of $\cal S^{-\nu}_\lambda(p,d-p)$ is non-empty, then the sign of $S(h)$ is constant and equal to the sign of $\lambda$.

\begin{figure}
\begin{tikzpicture}
\node at (0,0) {\includegraphics[width=14cm]{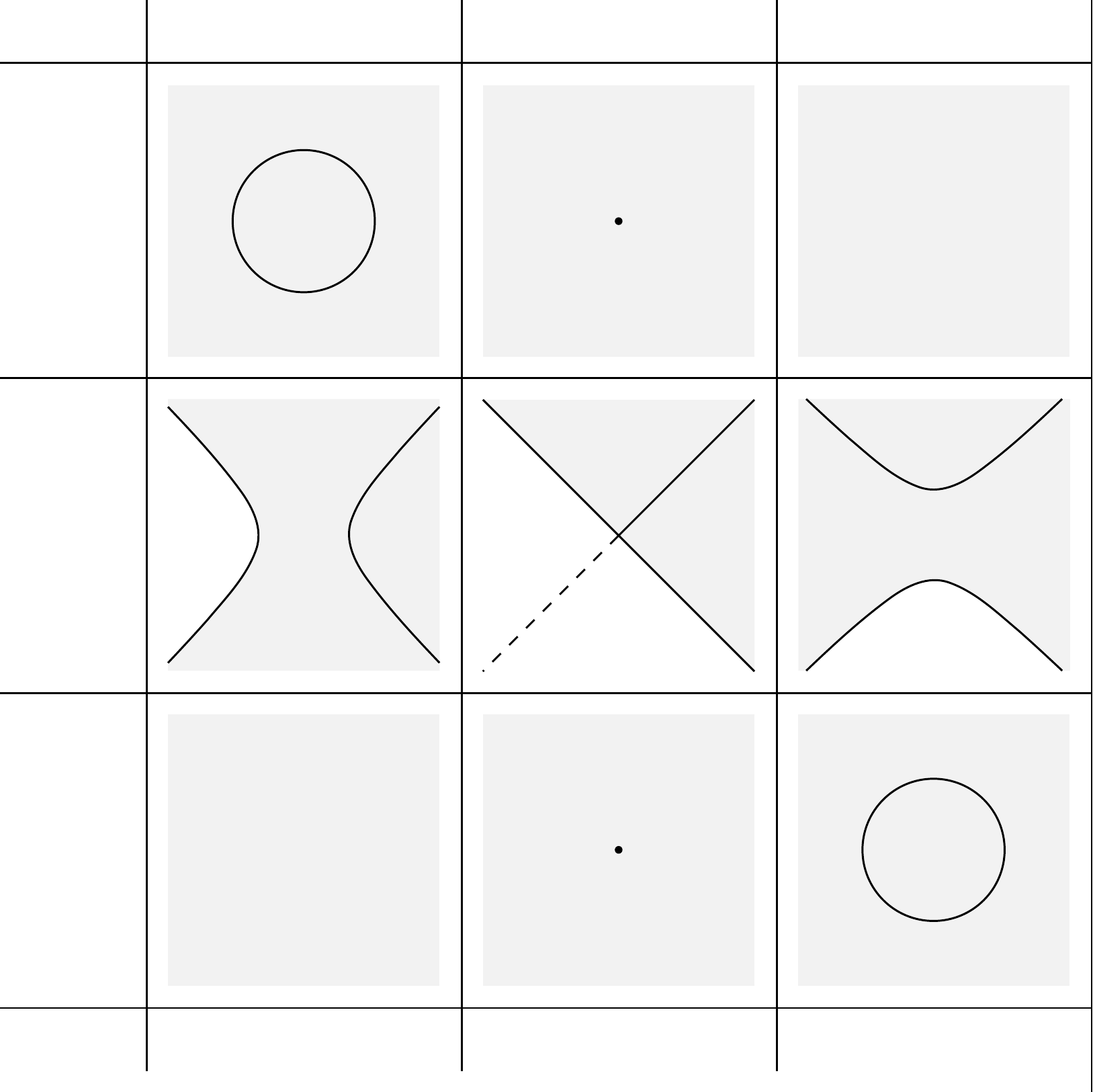}};

\node at (-3.1,6.6) {$\lambda = -1$};
\node at (0.9,6.6) {$\lambda = 0$};
\node at (4.9,6.6) {$\lambda = 1$};

\node at (-6.1,4.1) {$\sigma=(2,0)$};
\node at (-6.1,0.1) {$\sigma=(1,1)$};
\node at (-6.1,-3.9) {$\sigma=(0,2)$};

\node at (-3.1,-6.35) {$S(h)<0$};
\node at (0.9,-6.35) {$S(h)=0$};
\node at (4.9,-6.35) {$S(h)>0$};

\node at (-2.05,5.45) {$[+,-]$};
\node at (-3.1, 4.15) {$[-,-]$};
\node at (1.75,5.05) {$[+,0]$};
\node at (4.96,4.15) {$[+,+]$};

\node at (-3.1,1.4) {$[-,-]$};
\node at (-1.85,0.1) {$[+,-]$};
\node at (0.9,1.4) {$[-,0]$};
\node at (2,0.1) {$[+,0]$};
\node at (4.96,1.4) {$[-,+]$};
\node at (4.96,0.1) {$[+,+]$};

\node at (-3.05,-3.9) {$[-,-]$};
\node at (1.8,-3) {$[-,0]$};
\node at (4.96,-3.9) {$[+,+]$};
\node at (6.0, -2.6) {$[-,+]$};

\end{tikzpicture}
\label{fig:compare:S}
\caption{
Table of the different signs of the scalar curvatures $S(\nabla,g)$ of $\nabla$ relative to the flat pseudo-Riemannian metric $g$, and $S(h)$ of the isochrone pseudo-Riemannian metric $h$, for the 14 models $\cal S_\lambda^\nu(\sigma)$ in dimension $2$, indicated in each domain as $[\on{sign} S(\nabla,g),\on{sign} S(h)] = [\nu,\on{sign}\lambda]$.
See Section~\ref{sub:more:isoc} for an interpretation of these signs.
}
\end{figure}

In dimension $d=2$, we get the simpler expression
\[ S(h) = 16\lambda\psi^2
,\]
so the sign of the scalar curvature of $h$ is constant and the same as the sign of $\lambda$.
If moreover $\lambda=0$ then $S(h)=0$, which has the remarkable implication that $h$ is flat.
Actually, a flattening map, i.e.\ a (local) diffeomorphism  sending $h$ to $g$, can be explicitly determined:
\begin{itemize}
\item In signatures $(2,0)$ or $(0,2)$, identifying $E$ with $\C$, this is the non-injective map $z\mapsto 1/3z^{-3}$, whose expression in polar coordinates is $(r,\theta)\mapsto (1/3r^3,-3\theta)$.
\item In signature $(1,1)$ this is the diffeomorphism from $U$ to itself whose expression in polar hyperbolic coordinates is the same: $(r,\theta)\mapsto (1/3r^3,-3\theta)$. This map has a rational expression in terms of the coordinates $x_1,x_2$, which we omit here.
\end{itemize}
One may want to compare the sign of $S(h)$, which is the sign of $\lambda$ according to Remark~\ref{rem:ssng} with the sign of the relative scalar curvature $S(\nabla,g)$, which is $\nu$ according to Remark~\ref{rem:curvaturesign}.
Hence they do not necessarily match.
We found it interesting to present the comparison of the different possible signs in Figure~\ref{fig:compare:S}.
Note that they agree if and only if the domain contains $0$.

Let us repeat here the interpretation of the sign of scalar curvature in dimension $2$ that we gave in Section~\ref{sub:ScalarCurvature}.
Under the infinitesimal holonomy associated with small injective loops, vectors undergo an infinitesimal isometry. The latter displaces a vector $v$ along a direction which is tangent to the unit ``sphere'' $q(v)=\on{cst}$, in an orientation that depends on the orientation of the loop and on the product of the sign of $q(v)$ with the sign of the scalar curvature: they move in the same direction as the orientation of the loop if this product is positive, in the other direction otherwise.

\medskip

Returning to arbitrary dimension $d\geq 2$, denote $\psi_\lambda = q(x)+\lambda$.
If $\lambda \neq 0$, we have $h =  g /\psi_\lambda^4$ and $1/\psi_\lambda^4 = 1/(q(x)+\lambda)^4\sim 1/q(x)^4 = 1/\psi_0^4$ when $|q(x)|\to +\infty$.
Hence in some way, if the model is unbounded, the isochrone metric for $\lambda\neq 0$ is tangent at infinity to the isochrone metric for $\lambda=0$ (one can reach infinity by a path of finite length for $h$: it is enough to take a straight line whose direction $v$ satisfies $q(v)\neq 0$).
However one also observes that for $\lambda\neq 0$, $|S(h)(x)|\to +\infty$ when $|q(x)|\to +\infty$ whereas $S(h)=0$ for $\lambda=0$, so the tangency is not that good.

\begin{figure}
\begin{tikzpicture}
\node at (0,0) {\includegraphics[width=4cm]{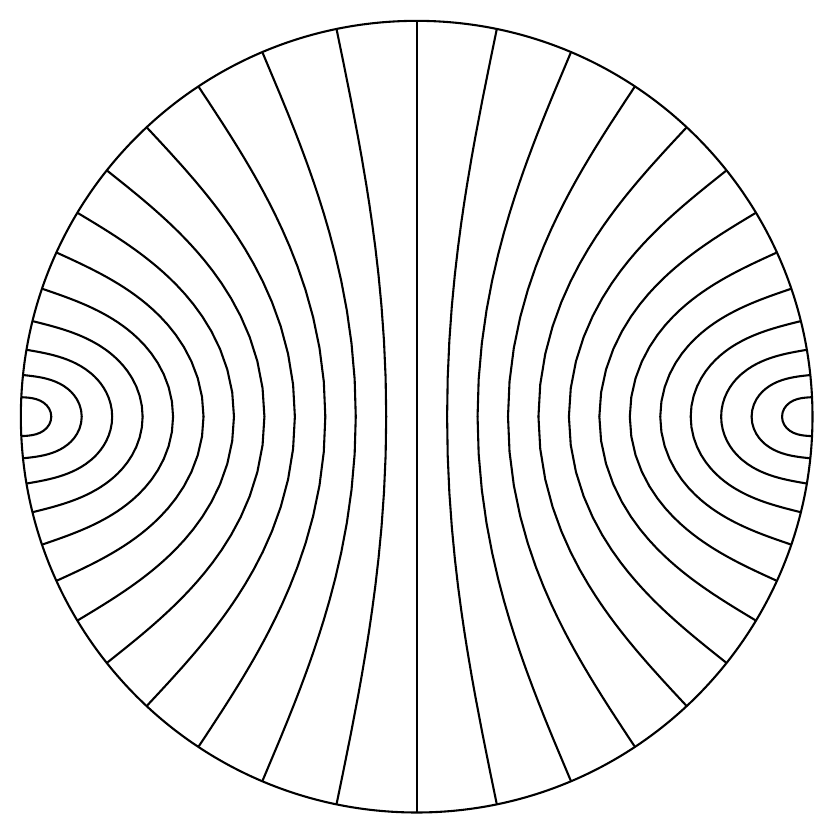}};
\node at (4.5,0) {\includegraphics[width=5cm]{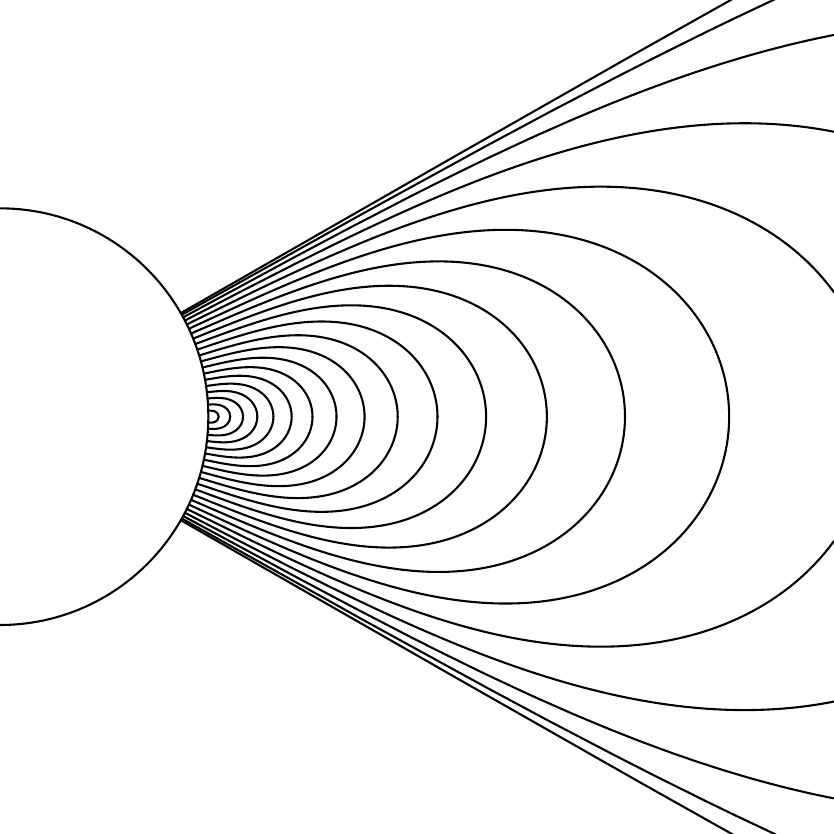}};
\node at (10,0) {\includegraphics[width=6cm]{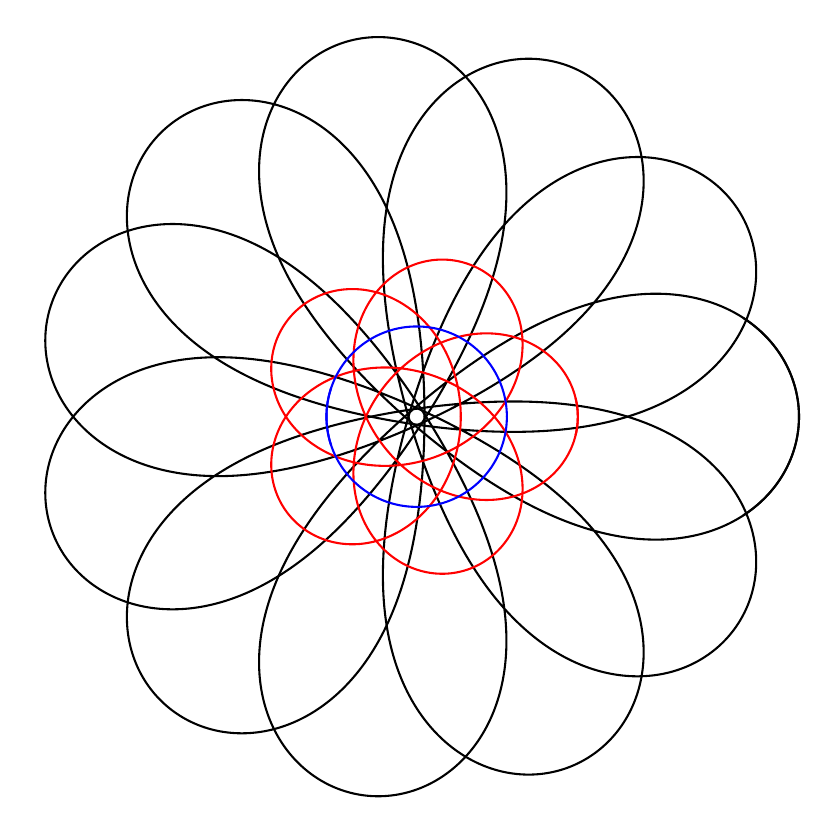}};
\end{tikzpicture}
\caption{A collection of geodesics of the isochrone metric for the models $\cal S^-_{-1}(2,0)$, $\cal S^+_{-1}(2,0)$ and $\cal S^+_{1}(2,0)$.
Middle: Only half of the unit disc is visible.
Left and middle: the fact that $S(h)\to 0$ at the boundary makes them quite straight near it.
Right: Most geodesics are quasiperiodic, filling a Euclidean annulus centered on $0$. We drew three periodic geodesics: in blue the only circular one (it has Euclidean radius $1/\sqrt{3}$), in red one with rotation symmetry of order 5, and in black one with rotation symmetry of order 11.}
\label{fig:geod:ex}
\end{figure}

By pure curiosity, we drew geodesics of the isochrone metrics of a few two-dimensional models in Figure~\ref{fig:geod:ex}.

\subsection{Forms and curvature}

We saw in Proposition~\ref{prop:Ricci} that the volume forms preserved by the connection of a canonical model are
\[
\frac{\beta}{\psi^{d+1}} dx_{1} \wedge \dots \wedge dx_{d}
\]
with $\beta \in \mathbb{R}^*$. We make here two remarks.

The canonical model $\mathcal{S}_{1}^{+}(d,0)$, i.e.\ the one whose domain is all of $\R^d$, has a finite total volume for any of these forms. We can check that for any other canonical model, the total volume is infinite.

The preserved volume forms $\omega$ are different from the volume form induced by the flat metric:
\[ \on{vol}_g = dx_1\wedge \cdots\wedge dx_d
\]
and from the volume forms induced by the isochrone metric $h=(1/\psi^4) g$:
\[ \on{vol}_h = \frac{1}{\psi^{2d}} dx_1\wedge \cdots\wedge dx_d.
\]

In dimension $d=2$, we would like to push a little bit further the analogy with the Levi-Civita connection of a Riemannian metric $\rho$.
More precisely, in dimension $2$ a particular case of the Gauss-Bonnet theorem can be stated in terms of a \emph{curvature form}
\[\kappa_\rho = G \on{vol}_\rho\]
where $G$ denotes the Gaussian curvature and $\on{vol}_\rho$ the volume form canonically associated with $\rho$: for any injective loop counterclockwise encompassing an oriented and simply connected domain $U$, the holonomy of the loop is a rotation whose angle, as measured per the orientation, is given by the integral of $\kappa_\rho$ on $U$ with respect to its orientation.

Such a strong statement cannot hold in our case since the local holonomy groups are the whole group $\on{SL}(2,\R)$ (Proposition~\ref{prop:LocHol}).
However, the first-order infinitesimal holonomies are rotations, and we computed in Proposition~\ref{prop:LocHol} that $M(1,2)$ has matrix $\frac{1}{\psi} J_\sigma$ with $J_\sigma = \begin{pmatrix} 0 & -\epsilon_2  \\ \epsilon_1 & 0 \end{pmatrix}$. Here $J_\sigma$ can be considered as a canonical infinitesimal rotation in the chart.
In view of the geometric interpretation of infinitesimal holonomies, it makes sense to define
\begin{equation}\label{eq:cf}
  \kappa = \kappa_\nabla = \frac{2\,dx_1\wedge dx_2}{\psi}
\end{equation}
and call it a \emph{curvature form} too, although it only satisfies Gauss-Bonnet at the first order (when the loop size tends to $0$).
We defined after Remark~\ref{rem:curvaturesign} the Gaussian curvature $G(\nabla,g)$ of $\nabla$ relative to $g$.
Note that $\kappa$ is also the product of this quantity with $\on{vol}_g$, i.e.\ $\kappa_\nabla = G(\nabla,g)\on{vol}_g$.
It would be tempting to define the \emph{relative curvature} as the quotient of the curvature form $\kappa$ by the invariant form $\omega$.
One problem is that the latter is only unique up to multiplication by a constant:
\[\omega = \frac{\beta}{\psi^3} dx_1\wedge dx_2\]
for any $\beta\neq 0$.
One immediately gets that the relative curvature is $2\psi^4/\beta$.
Remarkably, it vanishes on the boundary of the domain of the canonical model.

The curvature form of the isochrone metric is
\[\kappa_h = \frac{8\lambda}{\psi^2} dx_1\wedge dx_2.\]
It is thus different from the curvature form $\kappa_{\nabla}$ that we defined in Equation~\eqref{eq:cf}.

It is interesting to consider the following list:
\[
\begin{gathered}
\on{vol}_g = dx_1\wedge dx_2,\quad \kappa_\nabla=\frac{2}{\psi} dx_1\wedge dx_2,\quad \kappa_h=\frac{8\lambda}{\psi^2} dx_1\wedge dx_2
\\
\omega=\frac{\beta}{\psi^3} dx_1\wedge dx_2,\quad \on{vol}_h = \frac{\alpha^2}{\psi^4} dx_1\wedge dx_2
.
\end{gathered}
\]

\subsection{Parameterization of geodesics}\label{sub:geo:param}

Here we solve the geodesic equation for the hybrid connection $\nabla$ in canonical models. 
Recall that all geodesics of $\nabla$ are supported by straight lines $L$ of the pseudo-Euclidean space $E$ in which the model sits.
We will use the isochrone metric mentioned in previous sections.
We saw in Remark~\ref{rem:lcgeo} that for any Euclidean metric $\tilde g$ on $L$, the metric $\tilde g/\psi^4$ is isochrone, and this holds even if the direction of the geodesic lies in the light cone.
Choose any affine bijection $\phi:\R\to L$, $s\mapsto x_0 + se$, where $e\in\overrightarrow L$ is not zero.
Denote
\[ \psi_s=\psi\circ\phi
.\]
The index $s$ is purely symbolic and denotes the variable in which we work.
According to the point of view developed in Section~\ref{sub:isochrone}, $\psi_s$ is also a logarithmic potential of the connection $\nabla$ pulled back to $\R$ by $\phi$.
Geodesics supported by $L$ map by $\phi^{-1}$ to solutions $t\mapsto s(t)$ of
\[ \frac{\dot s}{\psi_s^2(s)} = \mathrm{cst}
\]
for some real constant depending on the initial speed of $s$ at time $t_0$.
In particular, any antiderivative $F$ of $1/\psi_s^2$ trivializes the connection, and sends the geodesics to linear maps (see also the proof of Proposition~\ref{prop:metr}). The same holds for any affine post-composition $aF+b$ of $F$.

The function $\psi_s$ develops as $\psi_s(s) = s^2 q(e) + 2(x_0\cdot e) s + q(x_0) +\lambda$, so is a polynomial of degree at most $2$ in $s$.
It cannot be the zero polynomial since it takes a non-zero value at the initial point.
The roots of $\psi_s$ are in correspondence with the intersections of $L$ with the equipotential $\psi=0$.

By a further affine change of variables $s\mapsto y$ we can reduce to studying the solutions of $\dot y = \mathrm{cst}/\psi_y(y)^2$, for the following five logarithmic potentials on $\R$:
\[
\text{(1)}\ \psi_y = 1 \quad
\text{(2)}\ \psi_y = y \quad
\text{(3)}\ \psi_y = y^2-1 \quad 
\text{(4)}\ \psi_y = y^2 \quad
\text{(5)}\ \psi_y = y^2+1
\]
This classifies geodesics into five types. One may subdivide (3) into further sub-cases: $|y|<1$ and $|y|>1$.

Denote $Q_0$ the equipotential $\psi=0$. 
Case~(1) occurs if and only if $q(e)=0$ and $(x_0\cdot e)=0$; equivalently iff.\ $L$ is directed by an isotropic vector and is disjoint from $Q_0$.
Case~(2) iff.\ $q(e)=0$ and $(x_0\cdot e)\neq 0$; equivalently iff.\ $L$ intersects $Q_0$ at a unique point counted with multiplicity, i.e.\ this point is different from the origin and the intersection is transverse.
Case~(3) iff.\ $L$ intersects $Q_0$ at two points.
Case~(4) iff.\ $L$ is tangent to $Q_0$ or if $\lambda=0$ and $L$ goes through $0$.
Case~(5) iff.\ $q(e)\neq 0$ and $L$ is disjoint from $Q_0$.

Trivializing maps $F$ are given by antiderivative of $1/\psi_y^2$ followed by any affine map. One can take:

\medskip

{\setlength{\extrarowheight}{4pt}%
\begin{tabular}{lll}
in Case (1):& $\psi_y = 1$,& $F=\on{id}$,
\\
in Case (2):& $\psi_y = y$,& $F(y)=-1/y$,
\\
in Case (3):& $\psi_y = y^2-1$,& $F(y) = F_{-1}(y) = -\frac{1}{y-1}-\frac{1}{y+1}+\ln\left|\frac{y+1}{y-1}\right|$, 
\\
in Case (4):& $\psi_y = y^2$,& $F(y) = F_0(y) = -1/y^3$, 
\\
in Case (5):& $\psi_y = y^2+1$,& $F(y) = F_1(y) = \frac{y}{1+y^2} + \arctan  y$.
\end{tabular}}

\medskip

\begin{figure}
\begin{center}
\includegraphics[scale=0.8]{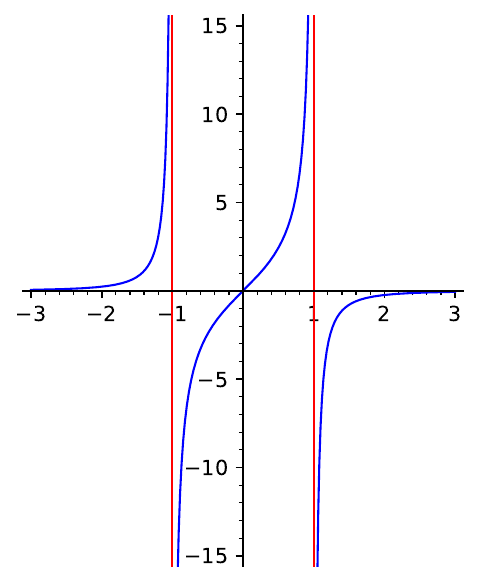}
\includegraphics[scale=0.8]{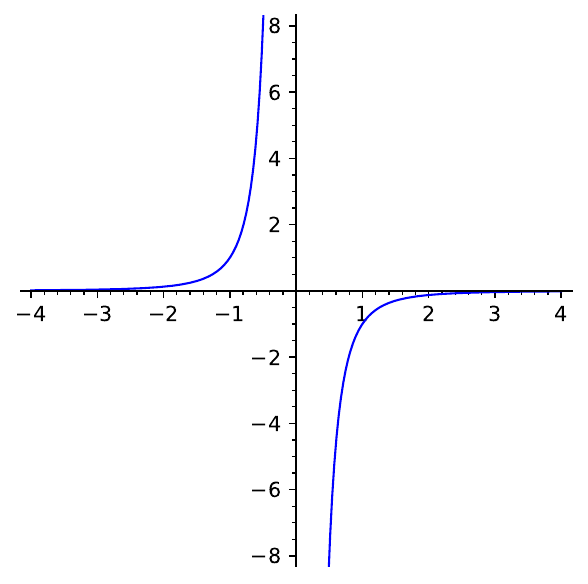}
\\
\includegraphics[scale=0.8]{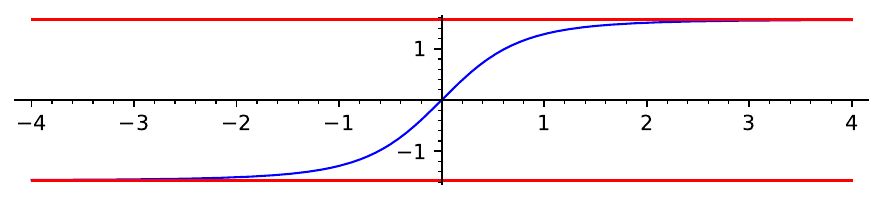}
\end{center}
\caption{Left: graph of the function $F_{-1}(r) = -\frac{1}{r-1}-\frac{1}{r+1}+\ln\left|\frac{r+1}{r-1}\right|$. Right: graph of the function $F_0(r) = -1/r^3$.
Bottom: graph of the function $F_1(r) = \frac{r}{1+r^2}+\arctan(r)$.}
\label{fig:ure}
\end{figure}

We chose the factor $a$ in each map $F$ so as to have a simple expression and positive derivative.
In each case, solving the geodesic equation amounts to solving
\[ F(y(t)) = \alpha t + \beta
,\]
i.e.\ to invert $F$, which is easy in Cases~(1),~(2) and~(4).
In Cases~(3) and~(5) this can be done numerically but seems hard or impossible to do symbolically.
The graphs of $F_{-1}$, $F_0$ and $F_1$ are shown in figure~\ref{fig:ure}.

Recall that $L = x_0 + \R e$ with $e\in \vec L$.
In Case~(1), $\gamma$ is an affine function of $t$.
In Case~(2) and~(4), the geodesic is defined on a semi-infinite interval of times. It tends in finite time to infinity and in infinite time to the unique intersection of $L$ with the equipotential $\psi=0$.
In Case~(5), the geodesic is defined on a bounded open interval of times, and tends to infinity at each end of this interval, so its support is the whole line $L$.
In Case~(3), there are two subcases: it is either defined for all times and tends in the future and in the past, respectively to the two intersections of $L$ with the equipotential $\psi=0$; or it behaves like in Case~(2).

In cases~(3),~(4) and~(5), when $y\longrightarrow \pm\infty$ we have $1/\psi_y^2\sim 1/y^4$.
It follows that $F(y)-F(\pm \infty)\sim \frac{-1}{3y^3}$ as $y\longrightarrow \pm \infty$.
In particular, any geodesic that tends to infinity along $L$ will do so in finite time with Euclidean distance satisfying $\sim \kappa/|t-t_0|^{1/3}$ for some $\kappa>0$.

\begin{rem}
For comparison, the geodesics in the Beltrami-Cayley-Klein model, and in the Beltrami-Poincaré model, also reach the boundary of the circle in infinite time.
However, in both cases, the Euclidean distance to the limit point is exponentially small with respect to the time, which is much faster than the decay rates $\kappa/|t|$ and $\kappa/|t|^{1/3}$ as $t\to\infty$ for geodesics of $\nabla$ tending to the boundary of a canonical model.    
\end{rem} 

\subsection{Example}\label{sub:ex:isoc}

Consider the model $\cal S^-_{-1}(2,0)$, defined on the unit disk $D$ of the classical Euclidean plane $\R^2$.
We saw that the geodesics of the associated connection $\nabla$, which are portions of straight lines in $\R^2$, have constant speed with respect to the isochrone metric $h = g/\psi^4$ where $g=dx_1^2+dx_2^2$ and $\psi(x) = x_1^2+x_2^2-1$.
Here we use Section~\ref{sub:geo:param} to detail the steps needed to compute the time $T$ it takes to go from one point $a\in D$ to another $b$ at speed $1$ with respect to $g$.

First, by definition
\[T=\int_{[a,b]} \frac{\|dx\|}{\psi(x)^2}.\] 
Consider the straight line $F$ in $\R^2$ passing through $a$ and $b$.
We are necessarily in case~3 of Section~\ref{sub:geo:param}.
The affine parameterization $y\in\R \to F$ leading to $\psi$ being proportional to $y^2-1$ is such that $-1$ and $1$ map to the intersections $u$ and $v$ of $F$ with $\partial D$.
Let $\ell = \|u-v\|/2$ be half the Euclidean length of $F\cap D$, so that $\|dx\|=\ell|dy|$.
Points $a$ and $b$ correspond to two parameters $y_a$ and $y_b$ in $(-1,1)$, and
\[ T=\int_{y_a}^{y_b} \frac{\ell |dy|}{\psi\left(\frac{(1-y)u+(y+1)v}{2}\right)^2}
.\]
Now $\psi\left(\frac{(1-y)u+(y+1)v}{2}\right) = \ell^2 (y^2-1)$ hence
\[T = \frac{|F_{-1}(y_a)-F_{-1}(y_b)|}{4\ell^3} 
\]
for the map $F_{-1}$ described before Case~3 of Section~\ref{sub:geo:param}.

It would be tempting to consider this time as a notion of distance between $a$ and $b$, however, it does not satisfy the triangle inequality as we show below with an example.

We consider $s\in(0,1)$ and the points $a=(s,0)$, $b=(0,s)$ and $o=(0,0)$.
Then
\[T(o,a)=T(o,b) = \frac{1}{4} F_{-1}(s),\qquad T(a,b) = \frac{1}{2\left(1-\frac{s^2}{2}\right)^{3/2}}F_{-1}(\frac{s}{\sqrt{2-s^2}}).\]
A numerical analysis shows that $T(a,b)\leq T(o,a)+T(o,b)$ if and only if $s\leq s_0 \approx 0.687$.
For instance, for $s=0.9$, we have $T(a,b) \approx 8.18$ whereas $T(o,a)+T(o,b) \approx 6.21$.
By decreasing the angle between $oa$ and $ob$, we can worsen the situation (this is consistent with the geodesics of $h$ being far from being straight chords of $D$, see Figure~\ref{fig:geod:ex}).

\subsection{Geodesically complete hybrid connections}\label{sub:geo:complete}

We remind the definition of geodesic completeness for arbitrary connections.

\begin{defn}\label{def:geo:complete}
A smooth manifold $M$ endowed with an affine connection $\nabla$ is \textit{geodesically complete} if every maximal geodesic of $\nabla$ is parameterized from $- \infty$ to $+\infty$.
\end{defn}

\begin{proof}[Proof of Theorem~\ref{thm:Uniqueness}]
Consider the Hessian manifold $(\mathcal{B}^{d},D,g)$ where $\mathcal{B}^{d}$ is the $d$-dimensional unit ball, $D$ is the standard affine connection and $g$ is the standard Euclidean metric.
Any Hessian potential for $g$ is of the form $\psi_{\lambda,A}:(x_{1},\dots,x_{d}) \mapsto \frac{1}{2}|x|^{2}+A(x)+\lambda$ for some $\lambda \in \mathbb{R}$ and some linear form $A$.
\par
Following Theorem~\ref{thm:MAIN}, any hybrid connection $\nabla$ for $(\mathcal{B}^{d},D,g)$ is determined by some hybrid potential $\psi = \psi_{\lambda,A}$.
Keeping in mind that $\frac{g}{\psi^{4}}$ is an isochrone metric for geodesics of $\nabla$ (see Theorem~\ref{thm:ThreeMetrics}), geodesics reach $\partial \mathcal{B}^{d}$ in infinite time only if $A=0$ and $\lambda=-1$.
It follows that the unique geodesically complete hybrid connection for $(\mathcal{B}^{d},D,g)$ is determined by the potential $\psi_{\lambda}:(x_{1},\dots,x_{d}) \mapsto \frac{1}{2}\left(|x|^{2}-1\right)$.
\end{proof}

\begin{rem}
The geodesics of the Cayley-Klein model are the usual straight lines and thus coincide with the (unparameterized) geodesics of the hybrid connection $\nabla$. However, the isochrone metric of $\mathcal{S}_{-1}^{-}(2,0)$ does not coincide with the Cayley-Klein hyperbolic metric nor with the Poincaré metric.
\end{rem}

\bibliographystyle{alpha}
\bibliography{bib}

\end{document}